\documentclass[aip,amsmath,amssymb,reprint,cha]{revtex4-1}
\usepackage{graphicx,float}
\usepackage{dcolumn}
\usepackage{bm}
\usepackage{xcolor}

\newcommand{\qed}{\hbox to 0pt{}\hfill$\rlap{$\sqcap$}\sqcup$\vspace{3mm}}
\usepackage{graphicx}

\newtheorem{thm}{Theorem}[section]
\newtheorem{theorem}[thm]{Theorem}
\newtheorem{lemma}[thm]{Lemma}

\newtheorem{proposition}[thm]{Proposition}
\newtheorem{assumption}[thm]{Assumption}
\newtheorem{remark}[thm]{Remark}
\newtheorem{example}[thm]{Example}
\begin{document}

\title[Stabilization of cycles with stochastic  control]{Stabilization of cycles with stochastic
prediction-based and target-oriented  control}


\author{E. Braverman}
\affiliation{Department of Mathematics and Statistics, University of Calgary, \\
2500 University Drive N.W., Calgary, AB T2N 1N4, Canada}
\email{maelena@ucalgary.ca}
\author{C. Kelly}
\affiliation{School of Mathematical Sciences, University College Cork,
Western Road, T12YT20, Cork, Ireland}
\email{conall.kelly@ucc.ie}
\author{A. Rodkina}
\affiliation{Department of Mathematics,   
The University of the West Indies, Mona Campus, Kingston 7, Jamaica}
\email{alexandra.rodkina@uwimona.edu.jm}



\date{\today}

\begin{abstract}
We stabilize a prescribed cycle or an equilibrium  of a difference equation using pulsed stochastic control. 
Our technique, inspired by Kolmogorov's Law of Large Numbers, 
activates a stabilizing effect of stochastic perturbation and allows for stabilization 
using a much wider range for the control parameter than would be possible in the absence of noise.  

Our main general result applies to both Prediction-Based and Target-Oriented Controls. 
This analysis is the first to make use of the stabilizing effects of noise for Prediction-Based Control; 
the stochastic version has previously been examined in the literature, 
but only the destabilizing effect of noise was demonstrated. 
A stochastic variant of Target-Oriented Control has never been considered, 
to the best of our knowledge, and we propose a specific form that uses a point equilibrium or one point on a cycle as a target. 
We illustrate our results numerically on the logistic, Ricker and Maynard Smith models from population biology.
\end{abstract}





\maketitle

\begin{quotation}
Various linear-type methods were developed to control otherwise unstable or chaotic behaviour of discrete maps.
Prediction-Based Control introduced by Ushio and Yamamoto in 1999 and two-parameter 
Target-Oriented Control proposed by Dattani et al in 2011 are among them. Stochastic perturbations 
were usually considered in two different contexts: as an intrinsic part of control which could 
diminish stabilization effects and thus should be kept in prescribed bounds, and as natural environmental noise which
may somehow control chaos. For instance, such noise can reduce oscillation amplitudes.

First, we consider control types incorporating both deterministic and stochastic components, both of which can 
have a stabilizing effect. Examples illustrate that stabilization can be achieved by noise but introduction of 
deterministic control, which cannot stabilize in itself, can influence the bounds for stabilizing stochastic 
perturbations. We analyze how the effective 
range of stabilizing control parameters may be extended by the introduction of noise. 

Second, stabilization of either an unstable equilibrium or an unstable orbit of a discrete equation is 
investigated. In addition, both regular (applied at each step) and pulsed (applied every $k$th step)
types of control are applied, and pulsed control can stabilize an equilibrium. 

Third, stochastic control is considered in very general settings. These results are later applied to 
Target-Oriented and Prediction-Based types of control involving stochastic component in the control parameter.
The application of stochastic Prediction-Based and Target-Oriented Controls for $k$-cycle stabilization is 
novel, stochastic versions of Target-Oriented Control have not been studied before.
\end{quotation}

\section[1]{Introduction}
\label{sec:int}

We investigate the use of pulsed stochastic control to stabilize a prescribed cycle of the difference equation 
\begin{equation}
 \label{eq:intr}
x_{n+1}=f(x_n),  \quad n\in {\mathbb N}_0, \,\,\quad x_0>0, 
\end{equation}
where  $\mathbb N_0:=\mathbb N\cup\{0\}$.
For a general class of control methods applied to \eqref{eq:intr}, 
we reduce this problem to the stabilization of a point equilibrium at zero and 
present a general theorem on pulsed stabilization of the zero equilibrium to equation \eqref{eq:intr}. 

We show how this theorem may be applied for two specific control methods: 
Prediction-Based Control (PBC) and a particular case of 
Target-Oriented Control (TOC). 
We introduce stochastic versions of both methods and study the interplay of the underlying control with 
stochastic perturbation of the control parameter, 
establishing results that show when the introduction of noise is beneficial for stabilization.  
In particular, we describe the stabilization of either a point equilibrium or a cycle by noise 
in the context of stochastic control. We also investigate the implications of pulsed control in this setting.

Our analysis allows us to demonstrate the stabilization of cycles for three commonly used models 
from population biology---Ricker, logistic and Maynard Smith---and we note that while stabilization 
of a point equilibrium by noise is quite a well-developed topic, stochastic stabilization of cycles is much less so.

Our technique is inspired by Kolmogorov's Law of Large Numbers which allows us to
 characterize the stabilizing effect of noise in our analysis, and may be stated as follows: 
\begin{lemma} \cite[page 391]{Shiryaev96} 
\label{thm:Kolm}
Let $(v_{n})_{n\in\ \mathbb N}$ be a sequence of independent identically 
distributed random variables where $\mathbb E |v_n|<\infty$, $n\in\mathbb{N}$. 
Denote the common mean $\mu:=\mathbb E v_n$, and the partial sum $S_n:= \sum_{k=1}^n v_k$.
Then $\lim_{n\to\infty}S_n/n=\mu$, a.s.
\end{lemma}

To the best of our knowledge, the idea of stabilization by noise goes back to 1950s for physical applications. 
Consider the well-known pendulum of Kapica~\cite{Kapica}, where stochastic perturbations can stabilize its top (otherwise unstable) position.  This stabilizing effect is stipulated by the type of noise and its intensity: noise that is too intense does not lead to stabilization, but noise that is insufficiently intense leads only to the preservation of the stability of the bottom equilibrium. For differential equations, a theoretical justification of stabilization by noise originated in the 1960s; see Hasminskii~\cite{Hasmin}. For both differential and difference equations, more detailed historical notes, as well as recent results on the topic are given in \cite{BKR2016,BR2019,Medvedev}. Recently, stability and stabilization of stochastic difference equations and systems, as well as cyclic and chaotic behaviour, has become a focus of many publications \cite{Bashk,Bashk1,Cortes,Guzik,Kadiev,Manjun,Rodina,Shaikhet}. Moreover, a developed theory of random difference equations was utilized to investigate differential equations \cite{Calatayud}, or discrete and continuous stochastic equations were considered in the framework of a single model \cite{Du}.

PBC was first introduced by Ushio and Yamamoto \cite{uy99} and was studied in detail in \cite{FL2010}.  
The case when the control is applied to \eqref{eq:intr} at every step, 
and the control parameter $\alpha\in (0,1)$ is subject to a stochastic perturbation, can be written
\begin{equation}
\label{eq:intro3}
x_{n+1}=  f(x_n)-\left(\alpha + l\xi_{n+1} \right) (f(x_n)-x_n),~ n\in {\mathbb N}_0,~x_0>0, 
\end{equation}
and was considered in \cite{BKR2016}. If the control is applied at every $k$th step, for $k \in \mathbb{N}$, $k>1$, 
then this is called {\em pulsed control}, and in the case of deterministic PBC it was investigated in \cite{BL2012,LizPotsche14}.

Previous results on PBC view stabilization as arising from the deterministic control and in spite of the presence of stochastic perturbations of low intensity. However, we can show that it is possible to stabilize the equilibrium of \eqref{eq:intr} by stochastic PBC \eqref{eq:intro3}, even for values of $\alpha$ which do not deliver stability in the absence of noise.  

TOC, applied to \eqref{eq:intr} at every step with target $T$, is characterized by
\begin{equation}
\label{TOC_eq}  
x_{n+1}=f\left( \alpha T+(1-\alpha)x_n \right), \quad T\ge 0, \quad \alpha \in [0,1).
\end{equation}
It was introduced in \cite{Dattani} and further investigated in \cite{Chaos2014,TPC}.   
In \cite{TPC} it was shown that TOC is topologically equivalent to the modified TOC equation
\begin{equation} 
\label{eq:intro5}
x_{n+1}=\alpha T+(1-\alpha)f(x_n), \quad T\ge 0, \quad  \alpha \in [0,1).
\end{equation} 

Note that in \eqref{eq:intro3} and other stochastic control models with $f: [0,\infty) \to [0,\infty)$, 
once the control $\alpha + l\xi_{n+1} \in [0,1]$, or, for TOC, the target is in addition non-negative, 
the expression in the right-hand side is non-negative. Assuming $ \alpha \in [0,1)$, $|\xi_{n+1}| \leq 1$, 
we get $l\leq \min\{ \alpha, 1-\alpha \}$. However, most of our results are local, 
and we can consider parameters outside of this domain, considering the truncated version when the right-hand side 
is a maximum of the computed value and zero, which is quite a typical approach in population ecology \cite{Schreiber2001}. 

The application of modified TOC to stabilize cycles in the context of higher order or vector difference equations was considered in 
\cite{BF2015,BF2017}. In the present article we find a relationship between the control parameter $\alpha$ and the noise intensity $l$, which guarantees local stability of a cycle after application of stochastic pulsed stabilization. 
In the case when the control parameter is such that the unperturbed model is stable, our method provides conditions on the noise intensity which preserve stability, similarly to \cite{BKR2016}. In the case when a deterministic system is unstable after application of the control, introduction of a noise with appropriate intensity guarantees stability. 

Pulsed control is essential in cases where application of control at each step is either impossible or inefficient from a practical or 
economical 
point of view. Here we consider linear types of control for nonlinear models with either one for PBC or two parameters for TOC involved.
Our control results are robust as a result of  the simplicity of the control structures and the continuity of the maps.
While stabilization of an unstable equilibrium with a control applied at every step is always possible, once the control intensity is sufficient, 
pulsed stabilization with PBC in the deterministic case is problematic \cite{BL2012}, sometimes it cannot be achieved for any values of 
the 
parameters. In our earlier paper \cite{BKR2016} a stochastic perturbation of deterministic PBC was explored but global stability was justified only 
in the case when all the values of the noisy control are within a range of parameters leading to stabilization in the deterministic case.
In contrast with these global results, highly local results are obtained in \cite{BR2019} where control and stabilization are achieved solely by 
noise but the neighbourhood of the equilibrium should be very small (sometimes less than $10^{-8}$) and thus anyway another control method is 
required in practical applications. Here our purpose is to combine a stabilizing effect of noise with PBC, expanding the range of allowable control 
parameters compared to \cite{BKR2016} and relaxing requirements on proximity of the initial value to the equilibrium in \cite{BR2019}.  

For TOC, delayed versions and stage-structed dynamics were considered in \cite{Chaos2014,BF2015,BF2017}, but
its stochastic version have never been explored. Here we obtain sufficient stabilization conditions for noisy TOC,
including its pulsed version. And again, noise can expand the range of parameters for which stabilization is achieved. From
this point of view, results of the present paper significantly generalize, for example, \cite{Chaos2014}.

All stochastic sequences considered in the paper are defined on 
a complete filtered probability space $(\Omega, {\mathcal{F}}$, $\{{\mathcal{F}}_n\}_{n \in \mathbb N}, {\mathbb P})$, 
where the filtration $(\mathcal{F}_n)_{n \in \mathbb{N}}$ is naturally generated by 
the sequence $(\xi_n)_{n\in\mathbb{N}}$, so that
$\mathcal{F}_{n} = \sigma \left\{\xi_{1},  \dots, \xi_{n}\right\}$.
We use the standard abbreviations ``a.s." for either ``almost sure" or ``almost surely" 
with respect to a fixed probability measure $\mathbb P$, and 
``i.i.d.'' for ``independent and identically distributed'', as it applies to sequences of random variables. 
A detailed discussion of relevant stochastic concepts and notation can be found, for example, in \cite{Shiryaev96}. 

Since equations in the present paper are motivated by population  models,  we will assume bounded stochastic perturbations in the following sense:

\begin{assumption}
\label{as:noise-1}
$(\xi_n)_{n\in \mathbb N}$ is a sequence of independent identically distributed  random variables,  
each satisfying $|\xi_n|\le 1$.
\end{assumption}

The paper is organized as follows. 
The main stabilization theorem is presented in Section ~\ref{sec:general} in its most general form.  A stochastic TOC method  is introduced and discussed in Section ~\ref{subsec:TOC}.  Results obtained in Section~\ref{subsec:PBC}  for stochastic PBC are generalizations of \cite{BKR2016}. In Section \ref{sec:exsim},  we illustrate some of our results with computer simulations. Section~\ref{sec:summary} contains a brief summary and discussion of potential directions for future research. All proofs are deferred to an Appendix in Section \ref{sec:ap}.


\section[2]{Local stabilization of a point equilibrium at zero by pulsed stochastic control}
\label{sec:general} 

In this section we present a generalized control theorem that will be applied to specific classes of model and control-type in the remainder of the article. Consider the difference equation 
\begin{equation}
 \label{eq:g1}
z_{n+1}=g(z_n),  \quad n\in \mathbb N_0, \,\,\quad z_0>0, 
\end{equation}
where the 
function $g$ satisfies a Lipschitz-type condition locally around zero:
 \begin{assumption}
\label{as:g}
For some $u_0>0$, there exists $L\ge 1$ such that
\begin{equation}
\label{cond:lipschg1}
|g(z)|\le L|z|, \quad |z|\le u_0.
\end{equation}
 \end{assumption}
Condition \eqref{cond:lipschg1} in Assumption \ref{as:g} is sufficient to ensure that equation \eqref{eq:g1} has a point equilibrium at zero, which we aim to stabilize by the application of pulsed stochastic control at each $k$th step, starting with the step $k-1$. In this article we are not concerned with the case where $L <1$, since it would immediately follow that $\lim_{n\to \infty}z_n=0$,  
for  $|z_0| \leq u_0$, and a control is unnecessary.
  
First, we characterize the control, which may depend on the function $g$, 
on a deterministic control parameter $\alpha\in[0,1)$, 
and on a coefficient $l>0$ describing the amplitude of a one-dimensional stochastic perturbation, 
satisfying Assumption \ref{as:noise-1}.  Note that if $\alpha=0$, any achievable control is due only to this perturbation.

Suppose we apply a general stochastic control to the right-hand side of \eqref{eq:g1}  at the $n=sk-1$ step, for each $s\in \mathbb N$, and represent the resulting stochastically controlled map by the  function $G: \mathbb R\times [0, 1]\times [0, l_0]\times [-1,1]\to \mathbb R$ for some $l_0>0$.  Then the stochastically controlled difference equation becomes (again with $z_0>0$)
\begin{equation}
\label {eq:genGalpha}
z_{n+1}=\left\{\begin{aligned} &G(z_n, \alpha, l, \xi_{n+1}),\quad n=sk-1, ~s \in {\mathbb N};\\
&g(z_n), \quad \mbox{otherwise}.
\end{aligned}
\right.
\end{equation}

Next, we place constraints on the form of the stochastically controlled map $G$ under which we will prove our main result in this section.

\begin{assumption}
\label{as:G}
Define the region $\mathcal B:=\{\alpha \in [0, 1), l\in [0, l_0], \, |v|\le 1\}$, and 
suppose that \eqref{cond:lipschg1} in the statement of Assumption~\ref{as:g} holds.
There exists a continuous function $\mathcal L: [0, 1]\times [0, l_0]\times [-1,1]\times [0, \infty)\to  (0, 
\infty)$ such that
\begin{enumerate}
\item [(i)] for $(\alpha, l, v)\in\mathcal B$, 
\begin{equation}
\label{cond:GLitsh}
\begin{array}{l}
|G(z, \alpha, l, v)|\le \mathcal L(\alpha, l, v, u)|z|, \\ |z|\le u\le u_0, ~~(\alpha, l, v) \in \mathcal B;
\end{array}
\end{equation}
\item [(ii)] for some $M>0$ \begin{equation}
\label{est:L}
\sup\{\mathcal L(\alpha, l, v, u): (\alpha, l, v)\in \mathcal B, \, |u|\le u_0\} = M;
\end{equation}
\item [(iii)] for $\xi$ satisfying Assumption \ref{as:noise-1}, for $L$ as given in \eqref{cond:lipschg1}, for $k$ as given in \eqref{eq:genGalpha}, and for some $\alpha \in [0, 1)$, $ l\in [0, l_0]$,
\begin{equation}
\label{cond:EGLitshu0}
\lambda:=-\mathbb E\ln \mathcal L(\alpha, l, \xi, 0)>(k-1)\ln L.
\end{equation}
\end{enumerate}
\end{assumption}

In our applications,  $G$ is the form which right hand side of the equation takes after the
shift of the equilibrium to zero, or after some other transformations,  and after application of
control.
The function $\mathcal L$ is a local Lipschitz constant of $G$ at zero.
Since the control is random, we have both $G$ and $\mathcal L$ random.

\begin{remark}
\label{rem:notauu}
Due to the continuity of  $\mathcal L$ in $u$, condition~\eqref {cond:EGLitshu0} 
implies that for some $u_1\in (0, u_0]$, 
$0\le u\le u_1$,
\begin{equation}
\label{cond:EGLitsh}
\lambda(u):=-\mathbb E\ln \mathcal L(\alpha, l, \xi, u)>(k-1)\ln L.  
\end{equation}
Inequality \eqref{cond:EGLitshu0} is the main assumption of the
paper,  and it immediately implies   \eqref{cond:EGLitsh}.
The fact that \eqref{cond:EGLitshu0}, \eqref{cond:EGLitsh} guarantee stability of the zero equilibrium 
is a consequence of Kolmogorov's  Law of Large 
Numbers. This applies to models that without stochasticity will be unstable, and can be connected to the illustration with Kapica's 
pendulum~\cite{Kapica}. Mathematically, this can be roughly described as possible decrease of Lyapunov exponents by introducing a random 
component with a 
zero mean. This approach goes back to H. Kesten in 1960-1970s \cite{FK,K}, see \cite{BR2019} for more details.
Note that conditions \eqref{cond:EGLitshu0}  and  \eqref{cond:EGLitsh} are quite close to necessary \cite{ABR2009}.
All the results of the present paper are proved under this condition, control
parameters for all the examples are chosen  to satisfy this assumption.
\end{remark}

Now we present the main result 
of 
this section.
\begin{theorem}
\label{thm:0equil}
Let Assumptions \ref{as:noise-1}, \ref{as:g}, \ref{as:G} hold, and let $\gamma\in (0, 1)$. 
Then there  exist $\delta_0>0$  and $\Omega_\gamma\subseteq \Omega $, $\mathbb P(\Omega_\gamma)>1-\gamma$, 
such that for each solution $z$ to Eq. \eqref{eq:genGalpha} with initial value $|z_0|\le \delta_0$ we have 
\[
\lim\limits_{n\to \infty}z_n(\omega)=0,~ \omega\in\Omega_\gamma,
\]   
where $z_n(\omega)$ is a sample path of the solution $z_n$. 
If we additionally suppose that $M L^{k-1}<1$ and $|z_0|\le\frac {u_0}{L^{k-1}}$, then 
\[
\lim\limits_{n\to \infty}z_n(\omega)=0,\quad \omega\in\Omega.
\]
\end{theorem}

\begin{remark}
\label{rem:nonrand}
The case  $M L^{k-1}<1$ can hold  when $G$ is nonrandom, or when $\mathcal L$ is nonrandom.  
In particular it is the case for PBC with a control that is either deterministic or subject to low-intensity stochastic perturbation, see \cite{BKR2016}. 
We discuss this in more detail in Sections \ref{subsec:TOC}, \ref{subsec:PBC}; see also Remarks \ref {rem:MLk<1} and \ref{rem:2PBC}.
\end{remark}

In this section we apply Theorem~\ref{thm:0equil} to stabilize cycles of
\begin{equation*}
x_{n+1}=f(x_n),  \quad n\in {\mathbb N}_0, \,\,\quad x_0>0.
\end{equation*}
First, we specify the structure of the map $f$. Then, we impose upon $f$ stochastic versions of 
TOC and PBC. The resulting stochastically controlled maps can be converted to form \eqref{eq:genGalpha}. 
For each model, we then derive assumptions on the control parameter $\alpha$ and the noise intensity $l$, which ensure condition~\eqref{cond:EGLitshu0} in the statement of Assumption \ref{as:G}, allowing us to apply Theorem \ref{thm:0equil}.

Suppose that $f$ is a real-valued and non-negative function possessing a cycle of period $d\in\mathbb{N}$, and it satisfies a Lipschitz-type condition locally around each point in the cycle:
\begin{assumption}
\label{as:fKd} 
For some $u_0>0$ and $d\in \mathbb N$, the continuous function $f: \mathbb R\to [0, \infty)$ is such that
\begin{enumerate}
\item[(a)]
$f(K_i)=K_{i+1}$ for  $i=1, \dots, d$, where $K_{d+1}:=K_1$;
\item[(b)] there exist $L_{i }>0$,  $i=1, \dots, d$, such that
\begin{equation}
 \label{cond:fLip}
 |f(x)-K_{i+1}|\le L_{i}|x-K_i|,~ 
x\in [K_i-u_0, K_i+u_0].
\end{equation}
\end{enumerate}
\end{assumption}

\begin{remark}
\label{rem:differentF}
If $f$ satisfies Assumption \ref{as:fKd}, it is not necessarily differentiable at $K_i$; consider for example $f(x)=|x|$. However, 
if $f'(K_i)$ exists, then  for each $u\in(0,u_0)$ and $x\in[K_i-u,K_i+u]$, \eqref{cond:fLip} is satisfied with
$L_i = |f'(K_i)|+\varepsilon (u)$, where  $\lim\limits_{u\to 0}\varepsilon(u)=0$.
\end{remark}

Notationally, set $f^2(x)=f(f(x))$, $f^j(x)=f(f^{j-1}(x))$, $j \in {\mathbb N}$ and note that under Assumption \ref{as:fKd}, each point of the set $\{K_1, K_2, \dots, K_d\}$ is an equilibrium for $f^d$. It follows that $f^d$ satisfies a generalized Lipschitz-type condition locally around each $K_i$:
\begin{lemma}
\label{lem:Lipschfd}
Let  Assumption \ref{as:fKd} hold and 
\begin{equation}
\label{def:L(d)}
L(d):=\prod_{i=1}^d \max\{1, \,L_i\}.
\end{equation}
Then,  for  $i=1, \dots d$,
\begin{equation}
\label{rel:Lipfd}
|f^d(x)-K_i|\le L(d)|x-K_i|  \quad \mbox{for} \quad |x-K_i|\le \frac{u_0}{L(d)}.
\end{equation}
\end{lemma}

Next we assume a Lipschitz-type relationship between the position of $x$ in the vicinity of a point in the $d$-cycle $K_i$, and the relative position of $f(x)$ to the next point in the $d$-cycle $K_{i+1}$.
\begin{assumption}
\label{as:fu0kcycle}
For some  $u_0>0$ and $d\in \mathbb N$ the function $f: \mathbb R\to [0, \infty)$  satisfies Part (a) of Assumption \ref{as:fKd}. There exist constants $\mathcal A_{i}\in \mathbb R$ and functions $\phi_i: \mathbb R\to \mathbb R$, $\psi_i: \mathbb R\to [0, \infty)$, $i=1, \dots, d$, such that 

\begin{enumerate}
\item [(i)] $\psi_i(u)\to 0$ as $u\to 0;$

\item  [(ii)] for each $u\in (0, u_0)$ and $x\in [K_i-u, K_i+u]$, $i=1, 2, \dots, d$,  
\begin{equation}
\label{cond:fudcycle}
\begin{split}
f^i(x)&=K_{i+1}+\mathcal A_{i}(x-K_i)+\phi_i(x);  \\
|\phi_i(x)|&\le \psi_i(|x-K_i|)|x-K_i|.
\end{split}
\end{equation}
\end{enumerate}
\end{assumption}
Since $K_{d+1}=K_1$ it follows that $x$ and $f^d(x)$ have a similar relationship in the vicinity of $K_1$:
\begin{lemma}
\label{lem:expfd}
Suppose that Assumption \ref {as:fu0kcycle} holds, and define for $u\in(0,1)$,
\begin{equation}
\label{def:A(d)u(d)}
\mathcal A(d):=\prod_{i=1}^{d}\mathcal A_{i},~ u(d):= \frac{u}{\prod_{i=1}^{d}\max\left\{ |\mathcal A_{i}|+\psi_i(u), \, 1 \right\}}.
\end{equation} 
Then,  there exist functions $\bar \phi: \mathbb R\to \mathbb R$ and   $\bar \psi: \mathbb R\to [0, \infty)$  such that for $x\in [K_1-u(d), K_1+u(d)]$, 
\begin{equation}
\label{cond:fudcycled}
\begin{split}
f^d(x)&=K_1+\mathcal A(d)(x-K_1)+\bar \phi(x),  \\
|\bar \phi(x)|&\le \bar \psi(|x-K_1|)|x-K_1|,
\end{split}
\end{equation} 
where $\bar \psi(u)\to 0$ as $u\to 0$.
 \end{lemma}
 

\subsection{Target-Oriented Control}
\label{subsec:TOC}

Deterministic modified TOC control is characterized in general by Eq.~\eqref{eq:intro5}. Consider a particular case when the target $T$ coincides with the equilibrium $K$ of $f$, and where the control parameter $\alpha$ is stochastically perturbed by an additive noise of intensity $l$. Then \eqref{eq:intro5} becomes
\begin{equation} 
\label{eq:intro7}
x_{n+1}=(1-\alpha-l\xi_{n+1})f(x_n)+(\alpha+l\xi_{n+1})K 
\end{equation}
for $n \in {\mathbb N}_0$. 
In fact, we apply \eqref{eq:intro5} not for all $n \in {\mathbb N}_0$ but at each $k$-th step, and we aim at either a point or a cycle stabilization. 
To the best of our knowledge, a combination of TOC with either pulsed control, stochastic control, or the use of part of a cycle as a 
target, is novel and we tackle here all three tasks.

In Sections \ref {subsubsec:TOCKk}  and \ref {subsubsec:TOCmd},  we present equations and conditions for the local stabilization of a 
point equilibrium and $d$-cycle, respectively, using stochastic TOC.  In Section \ref {subsubsec:TOCglob}  we investigate global stabilization.

\subsubsection{Pulsed stochastic TOC: stabilization of a point equilibrium.}
\label{subsubsec:TOCKk}

Suppose $f(K)=K$, for $K>0$. Consider the stochastic TOC model, pulsed at each $k$th step, with target $K$:
\begin{equation}
\begin{array}{l}
\label {eq:TOCKk}
x_{n+1}=\left\{ 
\begin{array}{l} (1-\alpha-l\xi_{n+1})f(x_n)\\ +(\alpha+l\xi_{n+1})K, \\   n=sk-1, ~~ s\in   \mathbb N,\\
f(x_n), \,\,  \mbox{otherwise},\\ 
\end{array}
\right.
\\
|x_0-K|<\delta. 
\end{array}
\end{equation}

If we denote
\begin{equation}
\label{def:subst1}
z_n:=x_n-K, \quad g(z):=f(z+K)-K,
\end{equation}
then \eqref{eq:TOCKk} takes the form of \eqref{eq:genGalpha} with 
\begin{equation}
 \label{def:GTOC}
 G(z_n, \alpha, l, \xi_{n+1})=(1-\alpha-l\xi_{n+1})g(z_n).
\end{equation}
Assuming that 
\begin{equation}
\label{cond:lipschgK}
|f(z+K)-K|\le L|z|
, \quad |z|\le u_0,
\end{equation}
we get  $\mathcal L(\alpha, l, v, u)= |1-\alpha-lv|L$, with $M=L\max_{|v|\le 1}|1-\alpha -lv|$, and condition~\eqref{cond:EGLitshu0} takes the form
\begin{equation}
 \label{cond:ELalkTOC}
 \lambda= -\mathbb E\ln |1-\alpha-l\xi|>k\ln L.
\end{equation}
Here we also assume that $\xi\neq\frac{1-\alpha}{l}$ in the case when  $\xi$ has a discrete distribution, see more details in 
\cite{BR2019}.
The following result now follows directly by an application of Theorem \ref{thm:0equil}:
\begin{theorem}
\label{thm:TOCKk}
Let Assumption \ref{as:noise-1}  and conditions~\eqref{cond:lipschgK}, \eqref{cond:ELalkTOC} hold.  
Then for each $\gamma\in (0, 1)$  
there exist $\delta_0>0$ and  $\Omega_\gamma\subset \Omega$ where $\mathbb P(\Omega_\gamma)>1-\gamma$,  such that  
for each solution  $(x_n)_{n\in\mathbb{N}}$ to equation~\eqref{eq:TOCKk} with initial value satisfying 
$|x_0-K|\le \delta_0$, we have  
\[
\lim\limits_{n\to \infty}x_n(\omega)=0, \quad \omega\in\Omega_\gamma.
\]   
 \end{theorem}
 
\begin{remark}
\label{rem:MLk<1}
Note that for  $\alpha<1-l$ we have $M=L(1-\alpha+l)$, and condition~\eqref{cond:ELalkTOC} holds if 
$1-\alpha+l\le L^{-k}$, which gives the following ranges for the parameter $\alpha$ and for the noise intensity $l$ to ensure stabilization:
\begin{equation}
\label{def:al}
\alpha\in \left( 1-L^{-k}, \, 1  \right), \quad  l\le \min \left\{ 1-\alpha, \, \alpha-1+L^{-k} \right\}.
\end{equation}
For a large value of $L$, the parameter $\alpha$ needs to be close to one, and  $l$ needs to be small. In this case, stabilization is due to the deterministic control $\alpha$, and only a small stochastic disturbance is allowed.

It is also possible to demonstrate active stabilization by noise, when $l$ is bigger than in \eqref{def:al}.
For Bernoulli distributed $\xi$ (taking each of the values $\pm1$ with the probability of 0.5), we have  $-\lambda=\frac 12\ln \left|(1-\alpha)^2-l^2\right|$,
and \eqref{cond:ELalkTOC} holds if 
\begin{equation}
\label{rel:lLTOC}
\sqrt{(1-\alpha)^2-L^{-2k}}<l<\sqrt{(1-\alpha)^2+L^{-2k}}.
\end{equation}
Note that \eqref{rel:lLTOC} can remain valid  even for large $L$ and  $\alpha=0$, but then $l<1$ should be close to 1.
\end{remark}

\subsubsection{Pulsed stochastic TOC: stabilization of a $d$-cycle.}
\label{subsubsec:TOCmd}
Let $k=md$ for some $m\in \mathbb N$, and let $f$ satisfy Assumption \ref{as:fKd}.  
Recall that each $K_i$, $i=1, \dots, d$, is a fixed point of $f^d$, and therefore of $f^{md}$. 
For simplicity we focus only on $K_1$, but our analysis applies equally to any other point in the cycle. 

Consider the equation
\begin{equation}
\label {eq:TOCmd}
x_{n+1}=\left\{\begin{aligned} &(1-\alpha-l\xi_{n+1}) f(x_n)+(\alpha+l\xi_{n+1}) K_1, \\ &  n=smd-1,~~  s\in   \mathbb N, \\
&f(x_n), \quad \mbox{otherwise},\quad   |x_0-K_1|\le \delta.
\end{aligned}
\right.
\end{equation}

Set $y_s:=x_{(s-1)md}$, $\bar \xi_s:= \xi_{(s-1)md}$, for $s\in\mathbb N$, and note that the sequence $(\bar\xi_s)_{s\in \mathbb N}$ 
satisfies Assumption \ref{as:noise-1}.  For  $ n=smd-1$,  $s\in\mathbb N$, we have 
\begin{equation}
\label{def:kmd}
\begin{array}{rl}
x_n=x_{smd-1}, & x_{n+1}=x_{smd}=y_{s+1}, ~~ y_1=x_0,\\
f(x_n) & = f^{md}(x_{n-md+1}) \\ & =f^{md}(x_{(s-1)md})=f^{md}(y_s). 
\end{array}
\end{equation}
Thus, \eqref {eq:TOCmd} can be transformed to 
\begin{equation}
\label{eq:1md}
\begin{array}{l}
\displaystyle
y_{s+1}=(1-\alpha-l\bar \xi_{s+1}) f^{md}(y_s)+(\alpha+l\bar \xi_{s+1}) K_1, \\ s\in\mathbb N, ~~  |y_1-K_1|\le \delta,
\end{array}
\end{equation}
which is in the form of \eqref {eq:TOCKk}  with  $k=1$, $K=K_1$,  $f^{md}$ instead of $f$, $y_1$ instead of $x_0$ and $s$ starting from 1.

Equation~\eqref{eq:1md}, in turn, can be transformed to \eqref{eq:genGalpha} if we set $$z_s:=y_s-K_1, \quad g(z):=f^{md}(z+K_1)-K_1.$$
Note that, by Lemma \ref{lem:Lipschfd}, $g$ satisfies \eqref{cond:lipschg1} with constant $L(md):=L^m(d)$ where $L(d)$ is defined by  \eqref {def:L(d)}.  Recall that $L(d)\ge 1$.  For $G$ defined as in \eqref{def:GTOC}, 
  $\mathcal L(\alpha, l, v, u)= |1-\alpha-lv|L^m(d)$, $M:=(1-\alpha +l)L^m(d)$, 
condition~\eqref{cond:EGLitshu0} 
takes the form 
\begin{equation}
 \label{cond:ELalkTOCd}
  \lambda:=-\mathbb E\ln |1-\alpha-l\xi|>m\ln L(d). 
\end{equation}
Therefore, Theorem \ref{thm:0equil}  implies $\lim\limits_{s\to \infty}y_s=\lim\limits_{s\to \infty}x_{(s-1)md}=K_1$, with any 
given probability $1-\gamma$ and small enough $\delta_0$. 

To extend this result to show that  $\lim_{n\to \infty}x_{nd+\bar j}=K_{\bar j}$, 
for each $\bar j=0, 1, \dots, d-1$, we require the next lemma:
\begin{lemma}
\label{lem:convmds+r}
Let $(x_n)_{n\in\mathbb{N}}$ be a solution of Eq.~\eqref{eq:TOCmd}. 
Let Assumptions \ref{as:noise-1}, \ref{as:fKd}, and condition~\eqref{cond:ELalkTOCd} hold.  
Then for each $\gamma\in (0, 1)$  there exist 
$\delta_0>0$ and  $\Omega_\gamma\subset \Omega$, with $\mathbb 
P(\Omega_\gamma)>1-\gamma$,  such that if $|x_0-K_1|\leq\delta_0$, and 
$s_0$ is such that $$|x_{smd}-K_1|<u_0L^{-m}(d) |1-\alpha+l|^{-m+1},\quad s\ge s_0,$$ then for $j=qd+\bar j$, $\bar j=0, 1, \dots, d-1$, $q=0, 1, \dots, m-1$, we have
\[
|x_{smd+j}-K_{\bar j+1}|\le |1-\alpha+l|^{m-1}L^{m}(d)|x_{smd}-K_{1}|,~ s\ge s_0.
\]
\end{lemma}
All the above brings us to the following theorem:
\begin{theorem}
\label{thm:TOCKd}
Let $(x_n)_{n\in\mathbb{N}}$ be a solution of Eq.~\eqref{eq:TOCmd}. 
Let Assumptions \ref{as:noise-1}, \ref{as:fKd}, and condition~\eqref{cond:ELalkTOCd} hold.  
Then for each $\gamma\in (0, 1)$  there exist $\delta_0>0$ and  $\Omega_\gamma\subset \Omega$, 
with $\mathbb P(\Omega_\gamma)>1-\gamma$,  such that if $|x_0-K_1|\leq\delta_0$,
\[
\lim\limits_{n\to \infty}x_{nd+\bar j}(\omega)=K_{\bar j},\quad\omega\in\Omega_\gamma,\quad \bar j=0, 1, \dots, d-1.
\]  
 \end{theorem}


\subsubsection{Global stabilization of a $d$-cycle by stochastic TOC}
\label {subsubsec:TOCglob}
Observe that if \eqref{cond:lipschgK} (when $d=1$, \eqref{cond:fLip} otherwise) 
holds globally on $\mathbb R$, condition~\eqref{cond:ELalkTOC} (when $d=1$, \eqref{cond:ELalkTOCd} otherwise) also holds. 
It is then possible to show that stochastic TOC \eqref{eq:TOCKk} (respectively \eqref{eq:TOCmd}) 
globally stabilizes the equilibrium $K$ (or a $d$-cycle).  
\begin{theorem}
\label{thm:TOCglob}
Theorem \ref {thm:TOCKk} (respectively Theorem \ref{thm:TOCKd}) holds for any $x_0>0$ if, in conditions~\eqref{cond:lipschgK}, \eqref{cond:ELalkTOC}  (respectively,  conditions~\eqref{cond:fLip}, \eqref{cond:ELalkTOCd}),  local Lipschitz constants are replaced with global Lipschitz constants.
 \end{theorem}
The proof modifies that of Theorem~\ref{thm:0equil} so that solutions are not required to stay in some neighbourhood of the initial 
value. Note however that the global Lipschitz constant $\bar L$ at the point $K$ (or in the case of a $d$-cycle, $\bar L_i$ at each $K_i$, $i=1,2 \dots, d$) can be quite large, reaching up to $\sup_{s\in \mathbb R}|f'(s)|$. Nonetheless, we will see in Example~\ref{ex:TOCRicker} that, in the case of Bernoulli $\xi$, condition~\eqref{cond:ELalkTOC} (respectively \eqref{cond:ELalkTOCd}) holds for large $\bar L$  (or  $\bar L^m(d)$) even with $\alpha =0$ if $l$ satisfies \eqref{rel:lLTOC}, where we replace $L$ by $\bar L$ (or in the case of a $d$-cycle, $L$ is replaced by $\bar L^m(d)$ and $k$ by $md$).

\subsection{Predictive Based Control}
\label{subsec:PBC}

The application of stochastic PBC is characterized in general by Eq.~\eqref{eq:intro3}. 
Following the order of investigation in Section \ref{subsec:TOC}, 
we will apply pulsed stochastic PBC at each $k$th step to stabilize a point equilibrium in Section \ref{subsubsec:PBCequil}, 
and to stabilize a $d$-cycle in Section \ref{subsubsec:PBCcycleI}.  

The results of this section are illustrated for a point equilibrium in Example \ref{ex:PBCBevH}, where even local stabilization is not possible for any $\alpha\in (0, 1)$ in the absence of a stochastic perturbation, and in Example \ref{ex:max}, where global stabilization is considered. The application of pulsed stochastic PBC to stabilize a 2-cycle is illustrated in Example~\ref{ex:PBCRicker}.

\subsubsection{Pulsed stochastic PBC: stabilization of a point equilibrium.}
\label{subsubsec:PBCequil}

Suppose $f(K)=K$ for $K>0$, and consider the stochastic PBC model, pulsed at each $k$th step:
\begin{equation}
\label {eq:kPBC}
x_{n+1}=\left\{\begin{aligned} &(1-\alpha-l\xi_{n+1})f(x_n)+(\alpha+l\xi_{n+1})x_{n}, 
\\ & n=sk-1, ~s\in   \mathbb N,\\
&f(x_n), ~~ \mbox{otherwise},~~ |x_0-K|\le \delta, 
\end{aligned}
\right.
\end{equation}
which, if we again use notation defined by \eqref{def:subst1}, takes the form of \eqref{eq:genGalpha}, with
\begin{equation}
\label{def:PBCG}
G(z, \alpha, l, v):=(1-\alpha-lv)g(z)+(\alpha+lv)z.
\end{equation}

We may identify constraints on $f$ that ensure condition~\eqref{cond:EGLitshu0} holds.
Suppose first that  for  $|x-K|\le u\le u_0$, $f$ admits expansion \eqref{cond:fudcycle} with $d=1$, $K_i\equiv K$,   $\mathcal A_i\equiv \mathcal A$,  $\phi_i(u)\equiv \phi(u)$,  and $\psi_i(u)\equiv \psi(u)$. 
In particular, this means that $f$ is differentiable at $K$ with derivative $\mathcal A$. Then
\begin{equation*}
\begin{split}
G(z_n, \alpha, l, \xi_{n+1}):=&\biggl[(1-\alpha-l\xi_{n+1})\mathcal A +\alpha+l\xi_{n+1}\biggr]z_n \\  &  +[1-\alpha-l\xi_{n}]\phi(z_n);\\
\mathcal L(\alpha, l, v, u):=&\left|\left(1- \alpha  \right) \mathcal A+\alpha +\left( 1-\mathcal A \right)lv\right|\\  &  +|1-\alpha+l|\psi(u);\\
 L:=&\mathcal A+\psi(u),\\
\end{split}
\end{equation*}
and conditions~\eqref{cond:EGLitshu0} and (by Remark \ref{rem:notauu})  \eqref{cond:EGLitsh} hold if 
\begin{equation}
\label{cond:mainalphau0}
-\mathbb E\ln \left|\left(1- \alpha  \right)\mathcal A+\alpha +\left( 1-\mathcal A\right)l\xi \right| >(k-1)\ln |\mathcal A|.
\end{equation}
The following theorem then follows immediately:
\begin{theorem}
\label{thm:mathcalA}
Let Assumption \ref{as:noise-1}, Assumption \ref{as:fu0kcycle} with $d=1$,    and condition~\eqref{cond:mainalphau0} hold, and let $(x_n)_{n\in\mathbb{N}}$ be a solution of \eqref{eq:kPBC}. 
Then for each $\gamma\in(0,1)$, there exist $\delta_0>0$ and $\Omega_\gamma\subset\Omega$ where $\mathbb{P}(\Omega_\gamma)>1-\gamma$, 
such that, if $|x_0-K|\le \delta_0$, we have 
\[
\lim_{n\to\infty} x_n(\omega)=0,\quad \omega\in\Omega_\gamma.
\]
\end{theorem}

\begin{remark}
\label{rem:mathcalA}
Relation \eqref {cond:mainalphau0} fails if $\mathcal A>1$ and  $l=0$ for any $\alpha\in (0, 1)$ and $k\in \mathbb N$.  
However, the presence of noise with $l>0$ can ensure local stability even for $\alpha=0$ and large $\mathcal A$.
To see this, assume that $\xi$ is  Bernoulli distributed. Then 
\[
\lambda = - \frac 12 \ln \left|\left[\left(1- \alpha  \right)\mathcal A+\alpha\right]^2 -\left( 1-\mathcal A\right)^2l^2 \right|, 
\]
and  \eqref {cond:mainalphau0} holds  if 
 $$
l_{low}:=\frac{\left[\left(1- \alpha  \right)\mathcal A+\alpha\right]^2 -\mathcal A^{-2(k-1)} }{\left( 1-\mathcal A\right)^2} <l^2
$$ $$< \frac{\left[\left(1- \alpha  \right)\mathcal A+\alpha\right]^2 + \mathcal A^{-2(k-1)}}{\left( 1-\mathcal A\right)^2}.
 $$ 
If  $\alpha=0$ and $k=1$, the lower bound on $l$ is given by $l_{low}=\frac{\mathcal A^2-1}{\left( 1-\mathcal A\right)^2}
 =1+\frac 2 {\mathcal A-1}>1$,  for each $\mathcal A>1$.   For  example, if $\mathcal A=2$, $\alpha=0$,  and $k=1$ we required $1.73\approx\sqrt{3}<l<\sqrt{5}$ for \eqref{cond:mainalphau0} to hold.  So it is reasonable to combine a nonzero control parameter $\alpha$ with nonzero noise intensity $l$.
 
Note that  for  $k=1$ and  any $\mathcal A\neq 1$ the lower bound on $l$ satisfies
$$
l_{low}
=\frac{\left(1- \alpha  \right)\left[\left(1- \alpha  \right)(\mathcal A-1)+2\right]}{\mathcal A-1}\to 0, 
\quad \mbox{as~~}  \alpha\to 1,
$$
while for $k>1$  
$$
l_{low}
= \frac{\left[\left(1- \alpha  \right)(\mathcal A-1)+1\right]^2-1}{\left( 1-\mathcal A\right)^2}+
\frac{1-\mathcal A^{-2(k-1)}}{\left( 1-\mathcal A\right)^2}
$$
$$
\to  \frac{1-\mathcal A^{-2(k-1)}}{\left( 1-\mathcal A\right)^2} ~~\mbox{as~~} \alpha\to 1.
$$
Therefore, when $k=1$ and  $\mathcal A\neq 1$,   for any $\varepsilon>0$  we can choose $\alpha\in 
(0, 1)$ such that the $l_{low}$ satisfies $l_{low}<\varepsilon$. In other words, small noise stabilizes the equilibrium if $\alpha$ is close to 1. When $k>1$ and $\mathcal A>2$, $$\frac{1-\mathcal A^{-2(k-1)}}{\left( 1-\mathcal A\right)^2}<\frac{1}{\left( 1-\mathcal A\right)^2}<1,$$ so there exists $\alpha\in(0,1)$ such that $l_{low}<1$.
\end{remark}

We can relax the assumption that $f$ is differentiable at $K$, instead requiring only that $f(x)-K$ changes sign from positive to negative as $x$ increases through some neighbourhood of $K$.  This corresponds to the case $\mathcal A<-1$ if $f$ is differentiable at $K$.  

\begin{theorem}
\label{thm:fK<>k}
Let $(x_n)_{n\in\mathbb{N}}$ be a solution to  \eqref{eq:kPBC} where \eqref{cond:lipschgK} holds and suppose that $f(x)>K$ for $x\in [K-u_0, K]$,  $f(x)<K$ for $x\in [K, K+u_0]$, and 
$ l\in (0, \min\{\alpha, \,1-\alpha\})$. 
\begin{enumerate}
\item [(i)]  If either
\begin{equation}\label{eq:tpia}
k=1,\quad \alpha>1-L^{-1},\quad l<L^{-1}-1+\alpha,
\end{equation}
or
\begin{equation}\label{eq:tpib}
\begin{split}
&k>1,~ 1<L^k<L+1,~ \alpha\in \left( 1-L^{-k}, L^{-k+1}\right),\\
&l\in \left(0, \min\left\{L^{-k}-1+\alpha,\quad L^{-k+1}-\alpha\right\}\right), 
\end{split}
\end{equation}
then for $|x_0-K|\le \delta\le u_0/L^k$,
\begin{equation*}
\lim\limits_{n\to \infty} x_n(\omega)=K,\quad\text{for all}\quad \omega\in\Omega.
\end{equation*}
\item [(ii)] If
\begin{equation}
\label {cond:max}
\begin{array}{ll}
\lambda:= & \displaystyle -\mathbb E\max\left\{\ln (|\alpha+l\xi|), \,  \ln \left( |1-\alpha-l\xi|L\right)\right\} \\
& >(k-1)\ln L,
\end{array}
\end{equation}
then, for any $\gamma \in (0, 1)$, there exist $\delta>0$ and $\Omega_\gamma\subset \Omega$ with $\mathbb P(\Omega_\gamma)>1-\gamma$, such that, 
for $|x_0-K|\le \delta$, we have 
$$
\lim\limits_{n\to \infty} x_n(\omega)=K,\quad \omega\in\Omega_\gamma.
$$
\end{enumerate}
\end{theorem}

\begin{remark}
\label{rem:2PBC}
In Theorem \ref{thm:fK<>k}, part (i) describes a situation where stabilization is due to the action of the underlying deterministic 
control, and the noise intensity $l$ is kept small to preserve this effect.  The case of  $k=1$  in part (i) of Theorem \ref{thm:fK<>k} 
where \eqref{eq:tpia} holds was covered in \cite{BKR2016}, and this analysis included establishing the global stability. However, for 
local stability Theorem  \ref{thm:fK<>k} generalizes the results of~~ \cite{BKR2016} to the case where the noise plays an active role 
in achieving stability: it applies to the situations where, for a chosen $\alpha\in(0,1)$ and $l=0$, the point equilibrium is unstable. The use of pulsed control here to achieve stabilization is also novel.

In part (ii) the noise also plays an active role, and \eqref{cond:max} gives a set of stabilizing parameters different from those in part (i). Here we present an example where \eqref{cond:max} is fulfilled, deferring a more detailed description and illustrative numerical simulation until Example \ref{ex:max}. 

The inequality  $\mathbb E\ln  [1-\alpha-l\xi]<k\ln L$ implies \eqref{cond:max} if 
$(1-\alpha-lv)L\ge \alpha+lv$ holds for each $|v|\le 1$. The latter is true when $
\alpha+l<1-(1+L)^{-1}$.  In the case of Bernoulli distributed $\xi$, this gives a lower bound for $l^2$ as $l_{\text{low}}=(1-\alpha)^2-L^{-2k}$.  It can be shown that the Ricker model with $r=2.41$, $L=1.5$, and the control with $\alpha=0.3$, $l=0.24$, $k=1$, satisfies \eqref{cond:max} in part (ii)  but not \eqref{eq:tpia} in part (i) of Theorem \ref{thm:fK<>k}.

Suppose more specifically that $(1-\alpha-lv)L= \alpha+lv$ for some $v\in(-1,1)$, and $\alpha+l>1-(L+1)^{-1}$. Then in order to satisfy  \eqref{cond:max} we need
$\alpha(1-\alpha)+l(1+l)<L^{-k}$.  It can be shown that for the Ricker model with $r=2.2$, the values $L=1.2$, $\alpha=0.28$, $l=0.27$, $k=2$ satisfy  \eqref{cond:max} in part (ii), but not \eqref{eq:tpib} in part (i) of Theorem  \ref{thm:fK<>k}.   More details may be found in Example~\ref{ex:max}.
\end{remark}

\begin{remark}
\label{rem:Schwarzian}
There are cases (applicable to both Ricker and logistic models) for which local stability implies global stability. 
Suppose $k=1$, so that control is applied at every step. Then, in the deterministic case, we have
\begin{align*}
f_{\alpha}(x)  & := (1-\alpha) x e^{r(1-x)} + \alpha x,\\ f_{\alpha}^{\prime} (x) & = (1-\alpha) (1-r x) e^{r(1-x)} + \alpha.
\end{align*}
The controlled map $f_{\alpha}$ is unimodal with a negative Schwarzian derivative, and so equilibria of the controlled deterministic equation are globally stable once they are locally stable. The general form of this result is due to Singer \cite{Singer}, see also \cite{Liz2007}. For deterministic PBC the result is in \cite{FL2010}, and some extensions of the idea can be found in \cite{Franco2020}.
The point equilibrium $K=1$ for the Ricker model is locally stable if
$$  
f_{\alpha}^{\prime} (1) = (1-\alpha) (1-r ) + \alpha > -1,
$$ 
or $\alpha \in (\alpha^*,1)$, where $\alpha^*=(r-2)/r$. 
According to \cite{BKR2016}, stabilization is achieved once $(\alpha - l,\alpha+l) \subseteq (\alpha^*,1)$.
\end{remark}

\subsubsection{Pulsed stochastic PBC: stabilization of a $d$-cycle.}
\label{subsubsec:PBCcycleI}
Suppose that $k=md$, $m\in \mathbb N$, and  Assumption \ref{as:fu0kcycle} holds. Consider the equation 
\begin{equation}
\label {eq:PBCmd}
x_{n+1}=\left\{\begin{array}{l} 
(1-\alpha-l\xi_{n+1})f(x_n) \\ +(\alpha+l\xi_{n+1}) 
 x_{n-md+1}, \\  n=smd-1,~  s\in   \mathbb N,\\
f(x_n), \quad \mbox{otherwise}. \\ 
\end{array}
\right.
|x_0-K_1|\le \delta.
\end{equation}
To this model, we apply transformation  \eqref{def:kmd} using notation as in Section \ref{subsubsec:TOCmd}, to get
$$
y_{s+1}=\left(1- \alpha - l\bar \xi_{s+1} \right)f^{md}(y_{s})+\left( \alpha + l\bar \xi_{s+1} \right) y_{s}, $$ $$  s\in   \mathbb N, \quad |y_1-K_1|\le \delta,
$$
which is covered by the case discussed in Section \ref{subsubsec:PBCequil}. To see this, substitute  $f^{dm}$ for $f$,  and $\bar \xi_s:= \xi_{sdm}$ for $\xi_n$ and $k=1$. Note that, by Lemma \ref{lem:expfd}, $f^{md}$ admits the expansion \eqref{cond:fudcycled}, substituting $md$ for $d$ and  $\mathcal A(md)=\mathcal A^m(d)$, where $\mathcal A(d)$ is defined as in \eqref{def:A(d)u(d)}. Therefore condition~\eqref{cond:EGLitshu0} has the form
\begin{equation}
\label{cond:ELalkPBCd} 
-\mathbb E\ln \left|\left(1- \alpha  \right)\mathcal A^m(d)+\alpha +\left( 1-\mathcal A^m(d)\right)l\xi \right| >0.
\end{equation}
Following the arguments of Section \ref{subsubsec:TOCmd}, we obtain the following theorem.
\begin{theorem}
\label{thm:PBCd}
Let $(x_n)_{n\in\mathbb{N}}$ be a solution of Eq. \eqref{eq:PBCmd}. Let Assumptions \ref{as:noise-1},  \ref{as:fu0kcycle}, and condition~\eqref{cond:ELalkPBCd}  hold.  Then for each $\gamma\in (0, 1)$  there exist $\delta_0>0$ and  $\Omega_\gamma\subset \Omega$, with $\mathbb P(\Omega_\gamma)>1-\gamma$,  such that if $|x_0-K_1|\leq\delta_0$,
\[
\lim\limits_{n\to \infty}x_{nd+\bar j}(\omega)=K_{\bar j},\quad\omega\in\Omega_\gamma,\quad \bar j=0, 1, \dots, d-1.
\]  
 \end{theorem}

\section[3]{Examples and Computer Simulations}
\label{sec:exsim}

In all simulations presented in this Section, we truncate the controlled map
$$ x_{n+1}=\max \left\{ (1-\alpha-l\xi_{n+1})f(x_n)+(\alpha+l\xi_{n+1})K, 0 \right\}  $$ 
in order to avoid negative values $x_n$.

We illustrate the results of Sections~\ref{subsec:TOC} and \ref{subsec:PBC} using difference equations associated with the Ricker function
\begin{equation}
\label{eq:ricker}
f_1(x) = xe^{r(1-x)}, \quad x\ge 0, 
\end{equation}
the logistic map
\begin{equation}
\label{eq:logistic}
f_2(x)=rx(1-x), \quad x\ge 0, 
\end{equation}
and the Maynard Smith model \cite{Thieme} with
\begin{equation}
\label{eq:BevHolt}
 f_3(x) = \frac{3x}{2+(x-3)^2}, \quad x\ge 0. 
\end{equation}
Here we simulate continuous uniformly distributed on $[0,1]$ and Bernoulli random variables $\xi_n$ to 
illustrate some cases from  Sections \ref{subsec:TOC} and \ref{subsec:PBC}. Each plot incorporating stochastic perturbations ($l> 0$) was generated with 3 runs, with a single run used to generate deterministic plots.

We start with the TOC method. Examples \ref{ex:TOCRicker} and \ref {ex:TOClogistic}
illustrate Theorem \ref {thm:PBCd}, $d=2$, applied to  Ricker and logistic  functions.

\begin{example}
\label{ex:TOCRicker}

First, we consider stochastic TOC \eqref{eq:TOCKk} applied at alternate steps ($k=2$) to a chaotic Ricker map $f_1$ 
satisfying \eqref{eq:ricker} with $r \approx 3.2716$ and a continuous uniformly distributed on $[-1,1]$ noise.
Without noise, $\alpha=0.7$ guarantees pulsed cycle stabilization, while uniformly distributed noise with $l=0.4$ leads to 
stabilization of $K=1$,
see Fig.~\ref{figure_revision}.

\begin{figure}[ht]
\centering
\includegraphics[height=.125\textheight]{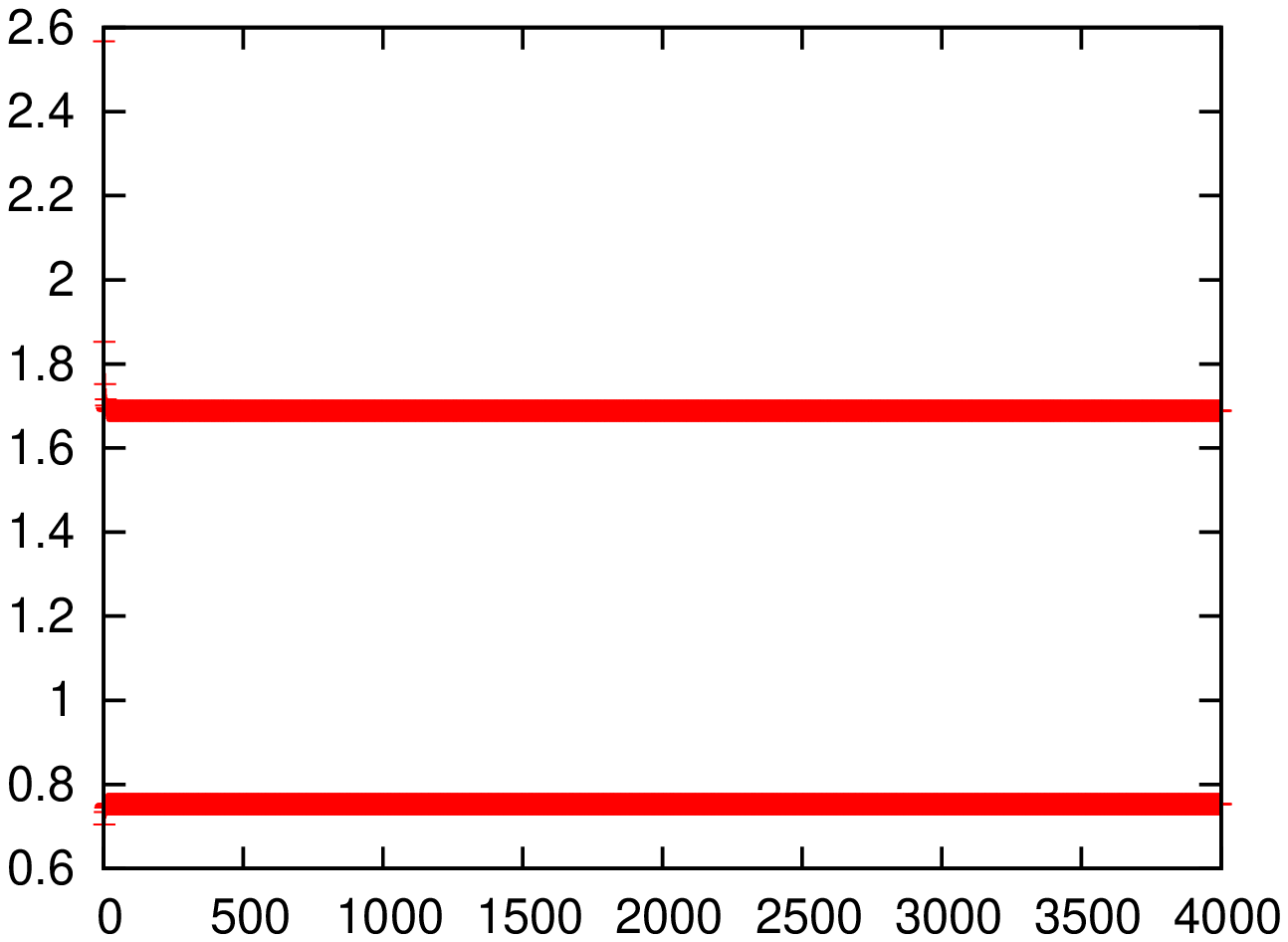}
\includegraphics[height=.125\textheight]{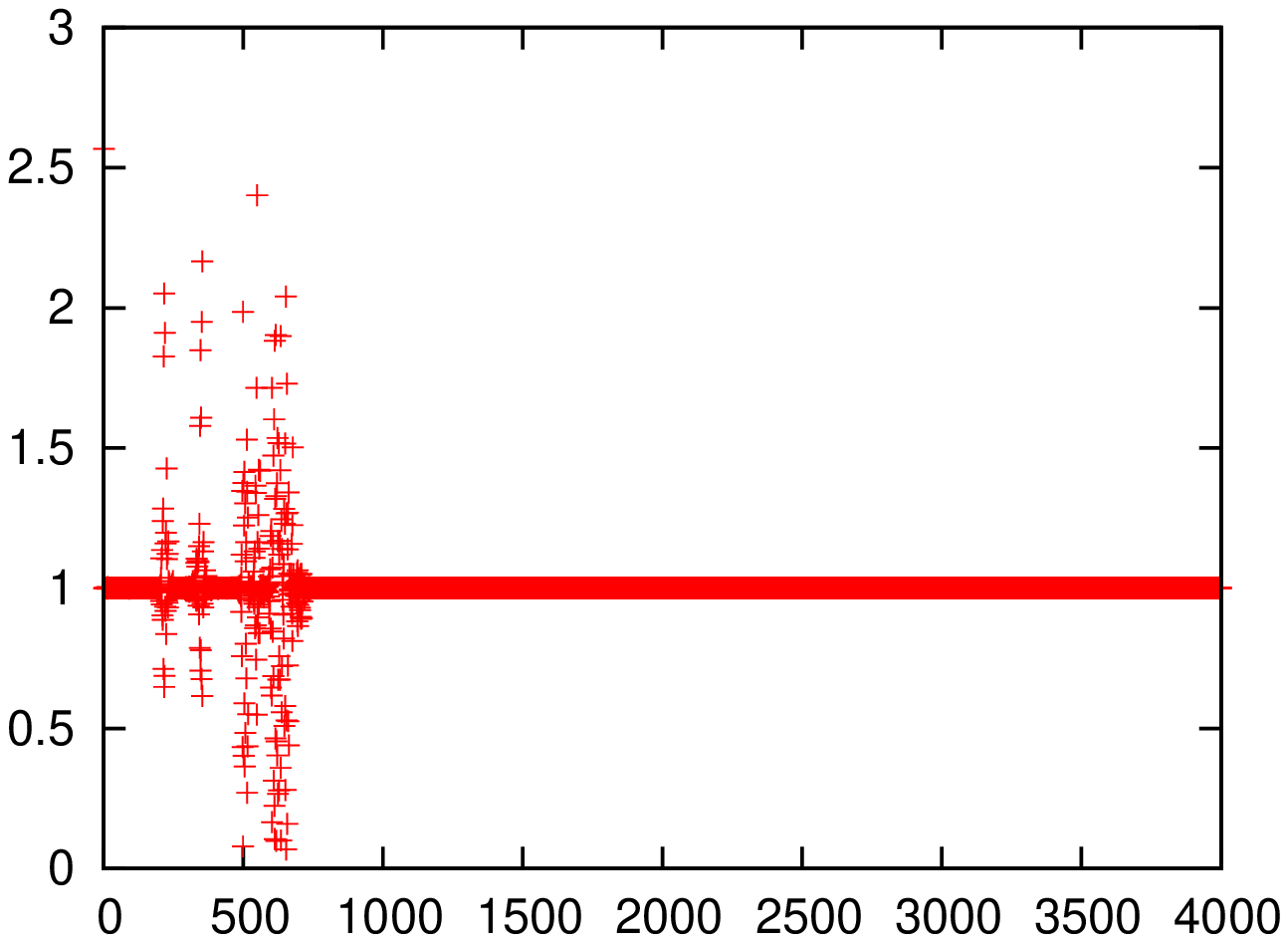}
\\
\includegraphics[height=.125\textheight]{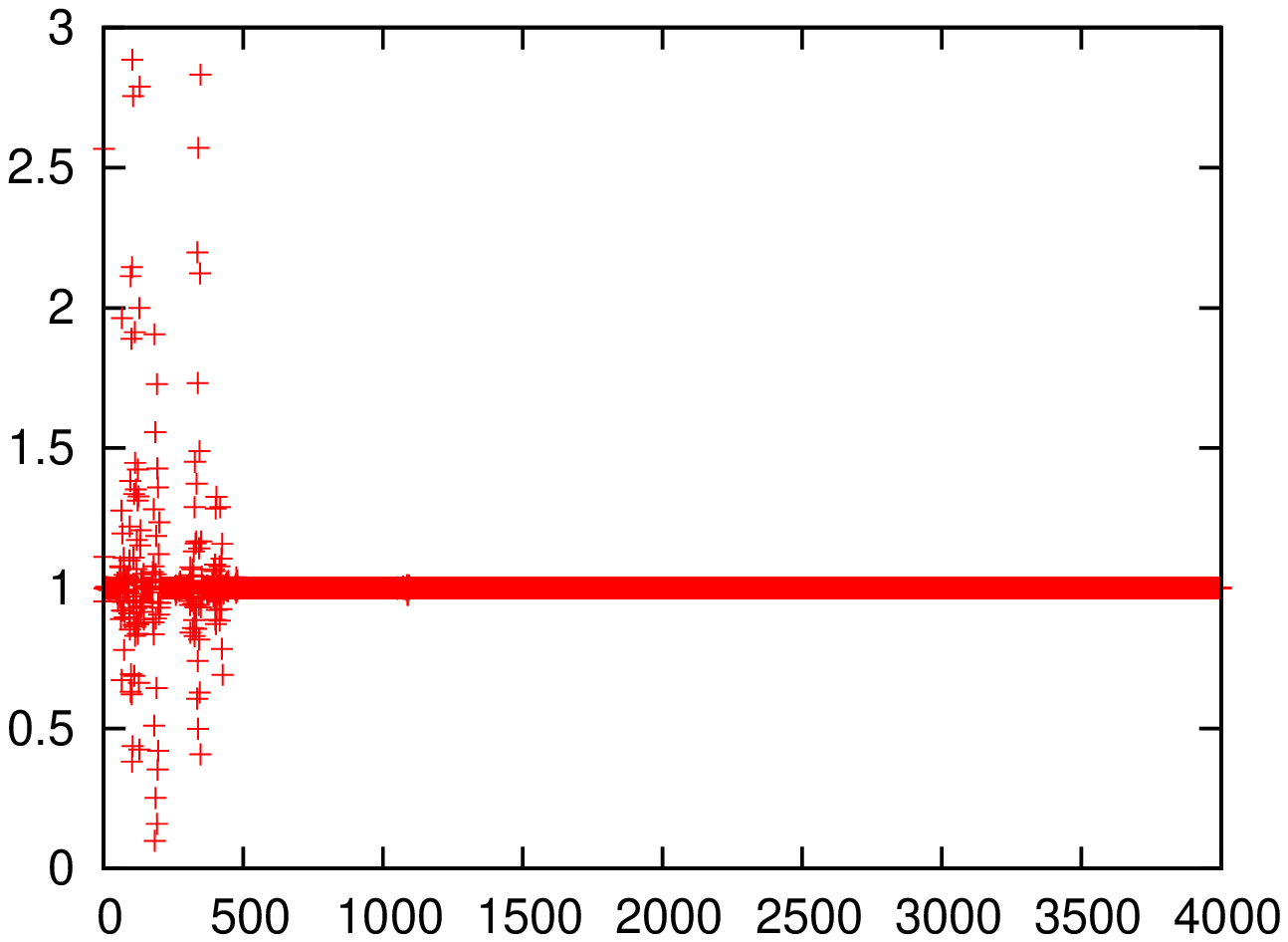}
\includegraphics[height=.125\textheight]{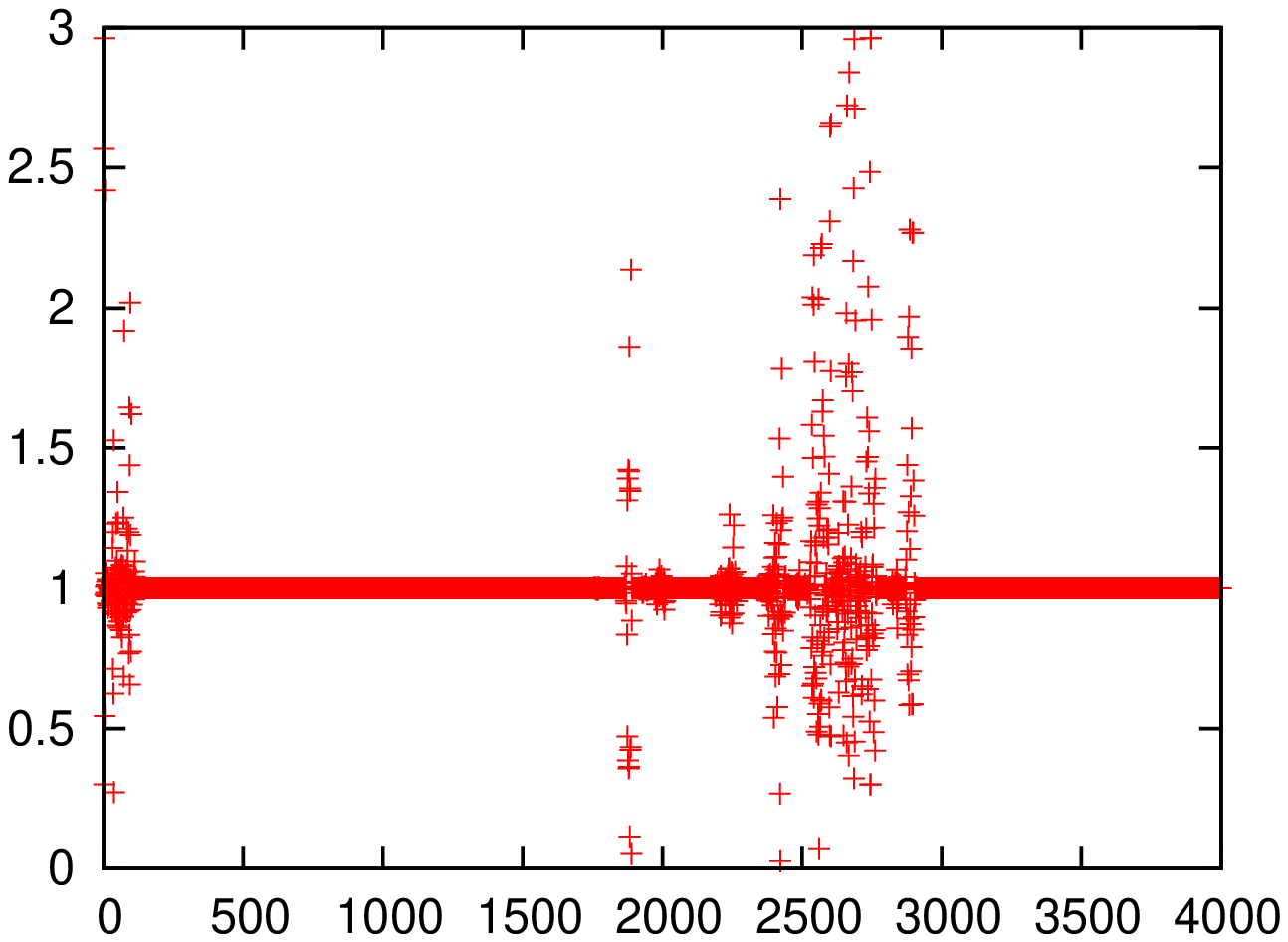}
\caption{  Model \eqref{eq:TOCKk} with  $f=f_1$ from \eqref{eq:ricker} with $r \approx 3.2716$, $\alpha=0.7$, $m=2$, $d=1$,
$x_0=0.5$ and (top, left) no noise, (top, right) $l=0.32$, (bottom, left) $l=0.4$ and (bottom, right) $l=0.42$.
}
\label{figure_revision}
\end{figure}

Next, apply a Bernoulli noise with a smaller $\alpha=0.3$ to stabilize 
2-cycle $K_1=0.1$, $K_2=1.9$ using $K_1$ as the target.
Fig.~\ref {figure6new}, right, presents stabilization for noise intensity $l=0.7$. 
For $l=0.7$, Fig.~\ref {figure6new}, left illustrates that there is no convergence to this 2-cycle, 
and similar results are obtained for $l>0.75$,
the range of values of $l$ that allow  stabilization is quite narrow.
\end{example}

\begin{figure}[ht]
\centering
\includegraphics[height=.125\textheight]{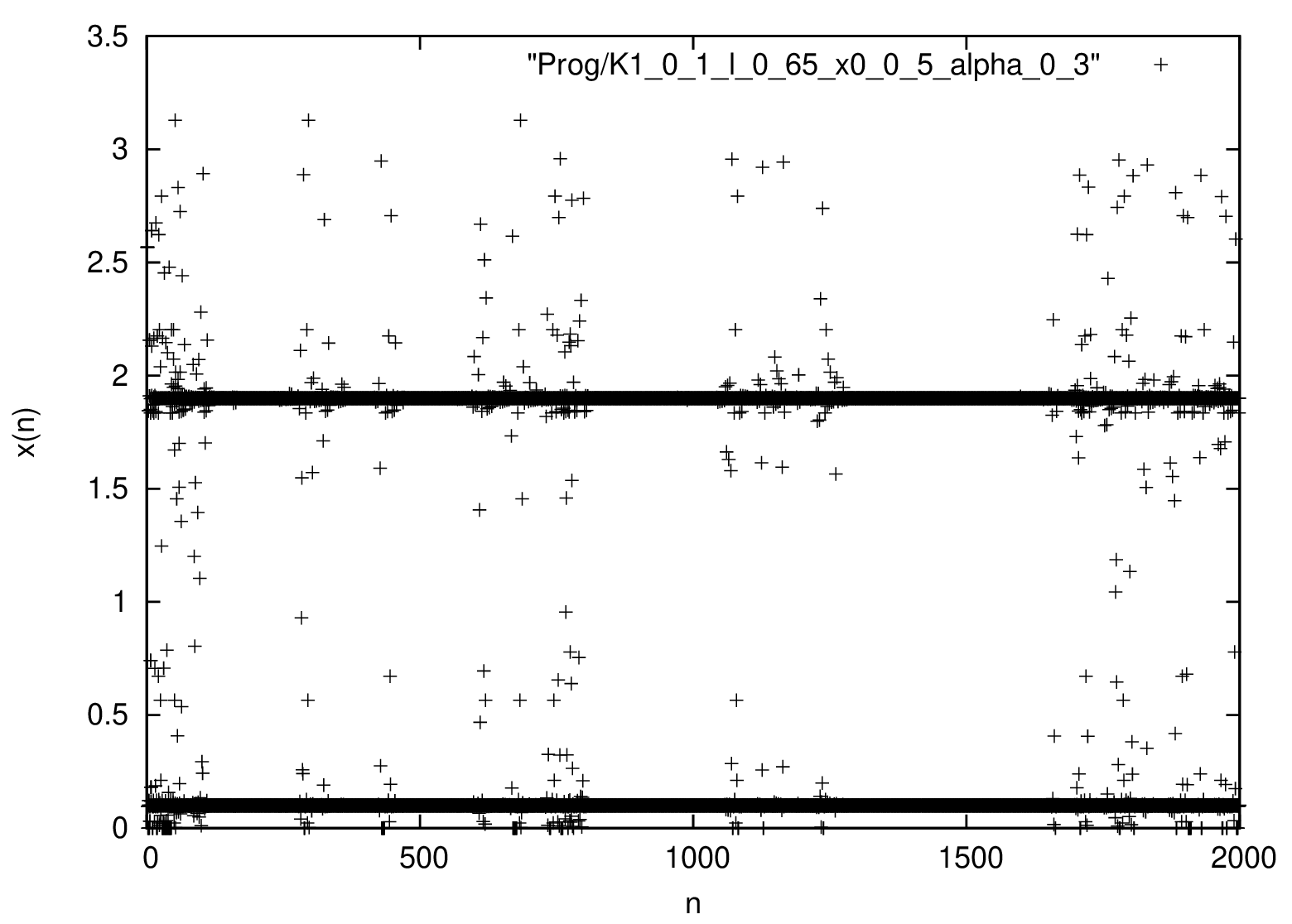}
\includegraphics[height=.125\textheight]{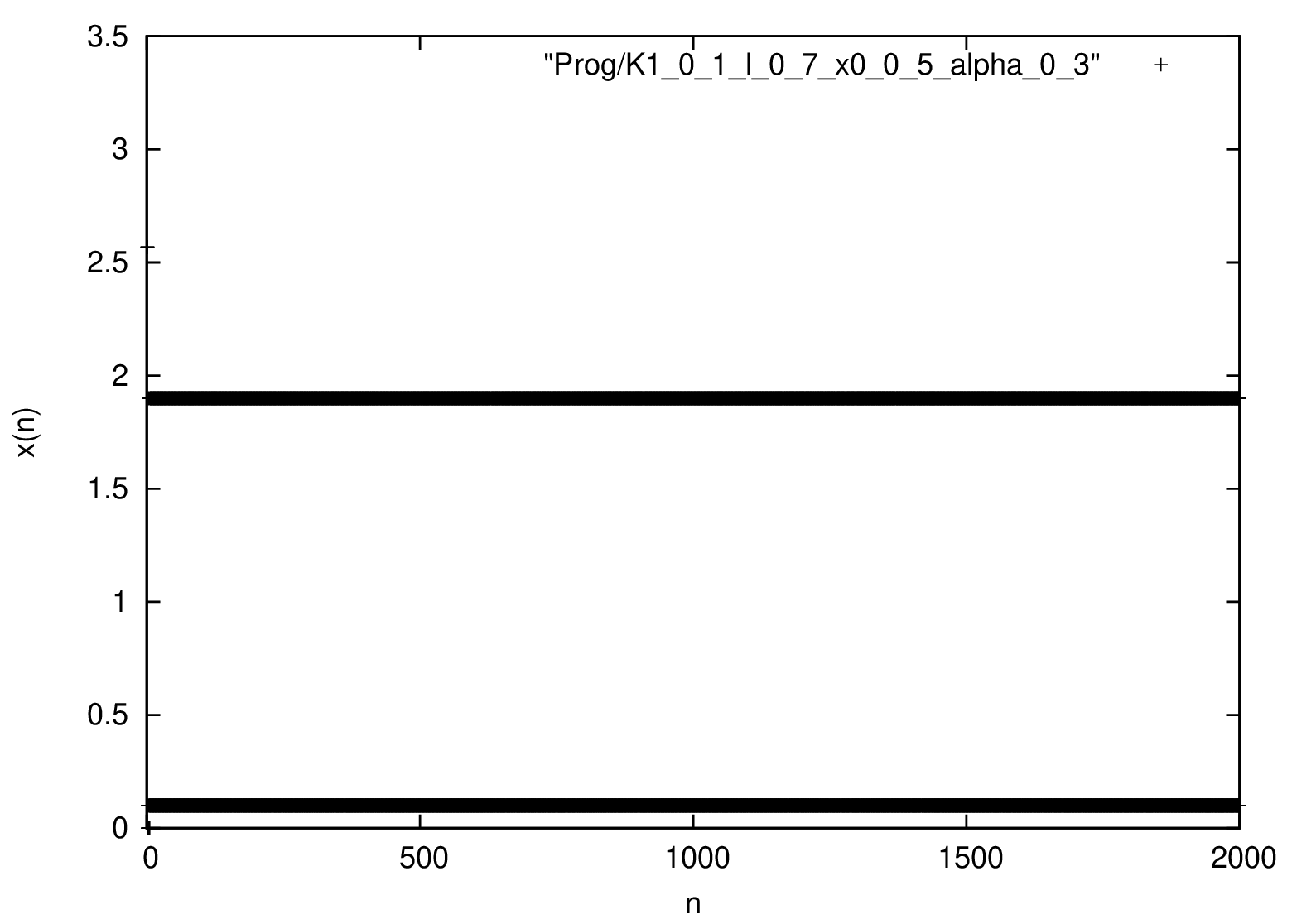}
\caption{Model (\protect{\ref{eq:TOCmd}}) with $f=f_1$ from \eqref{eq:ricker}  
with $r \approx 3.2716$, $\alpha=0.3$, $m=1$, $d=2$, $x_0=0.5$,  and
$l=0.65,0.7$. 
}
\label{figure6new}
\end{figure}

\begin{figure}[ht]
\centering
\includegraphics[height=.125\textheight]{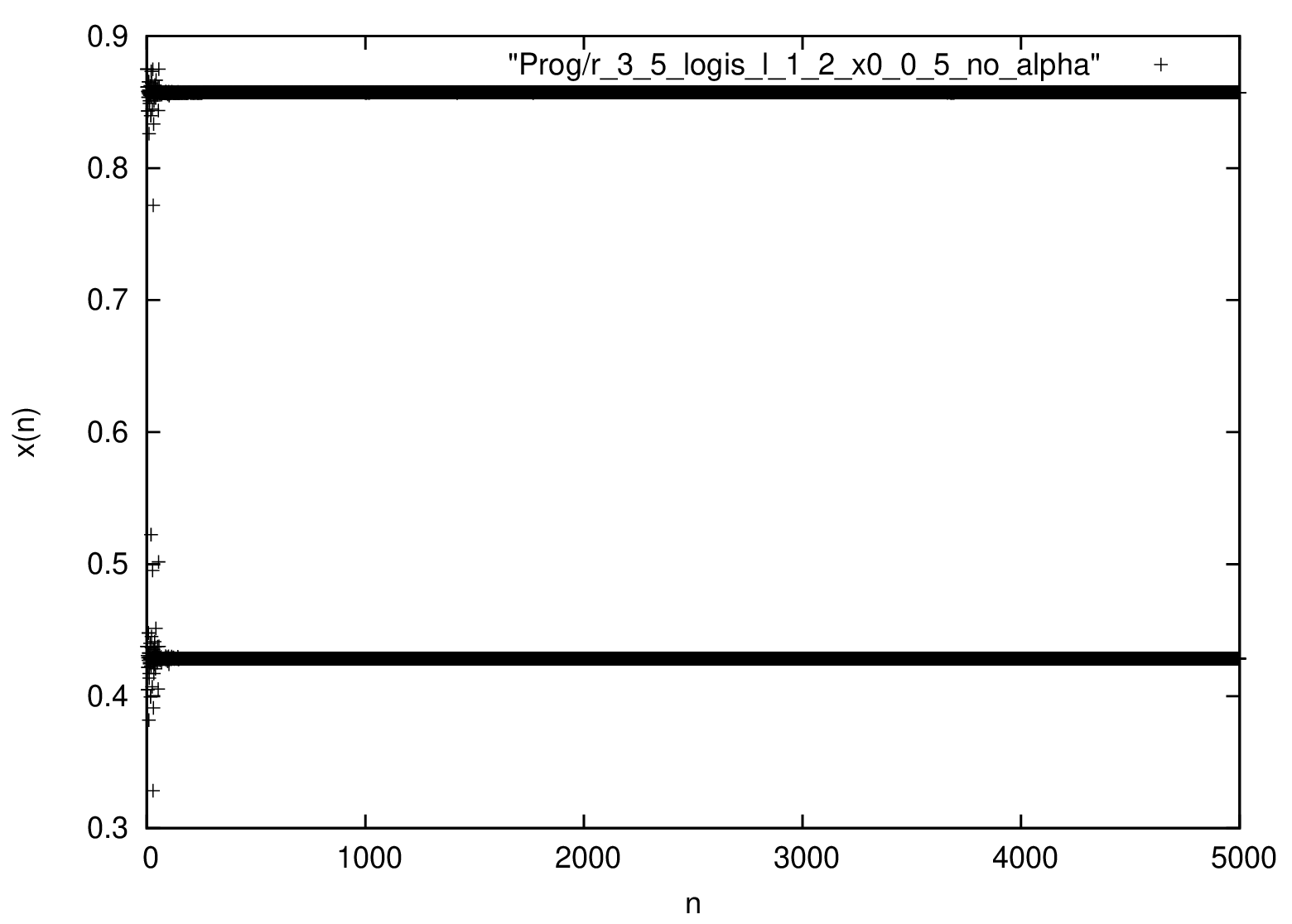}
\includegraphics[height=.125\textheight]{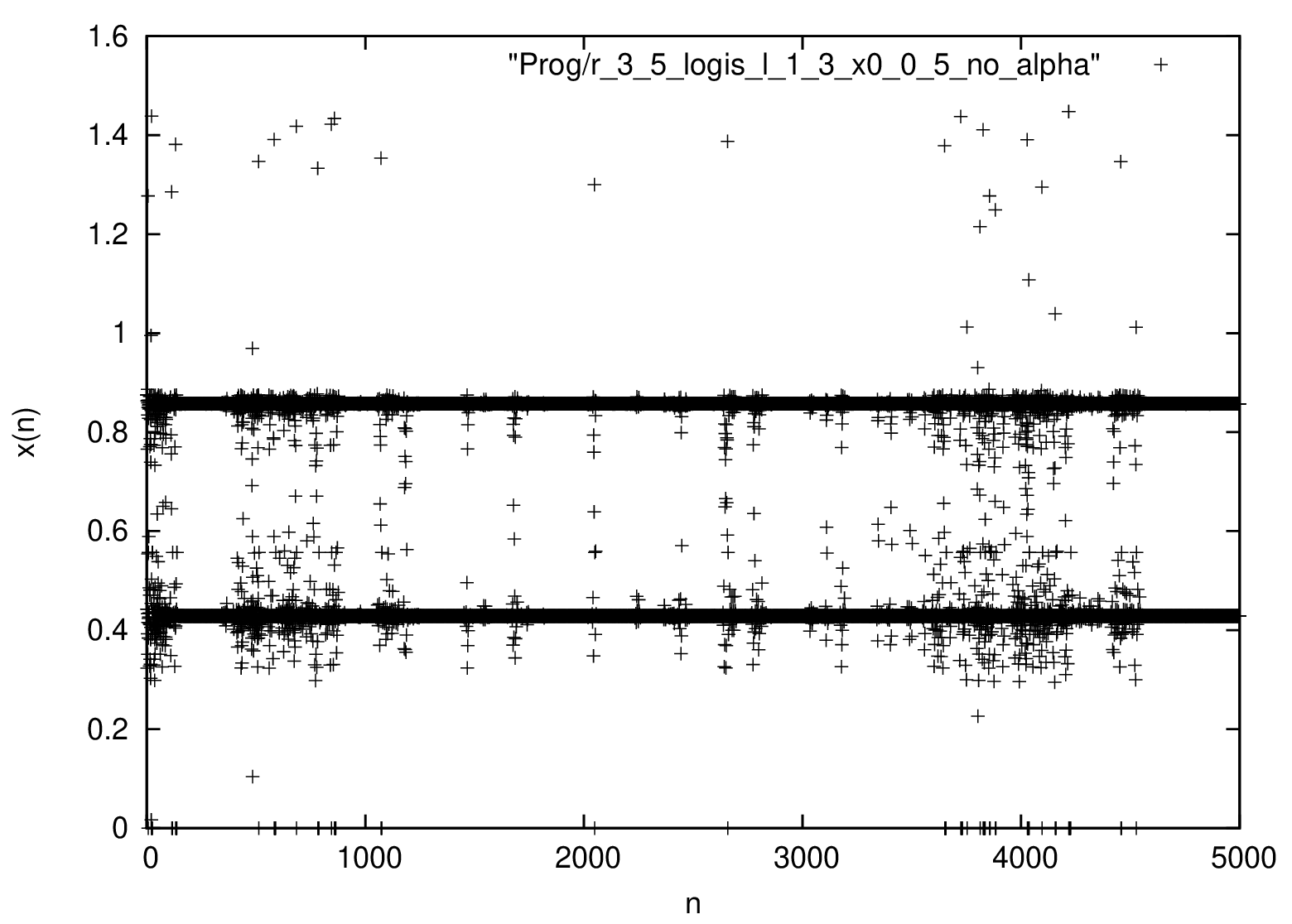}
\caption{
Model \eqref{eq:TOCKk} with $f=f_2$ from \eqref{eq:logistic},  $r = 3.5$, $\alpha=0$, $m=1$, $d=2$, $x_0=0.5$, and (from left to 
right) $l=1.2, 1.3$.
}
\label{figure7new}
\end{figure}

\begin{figure}[ht]
\centering
\includegraphics[height=.125\textheight]{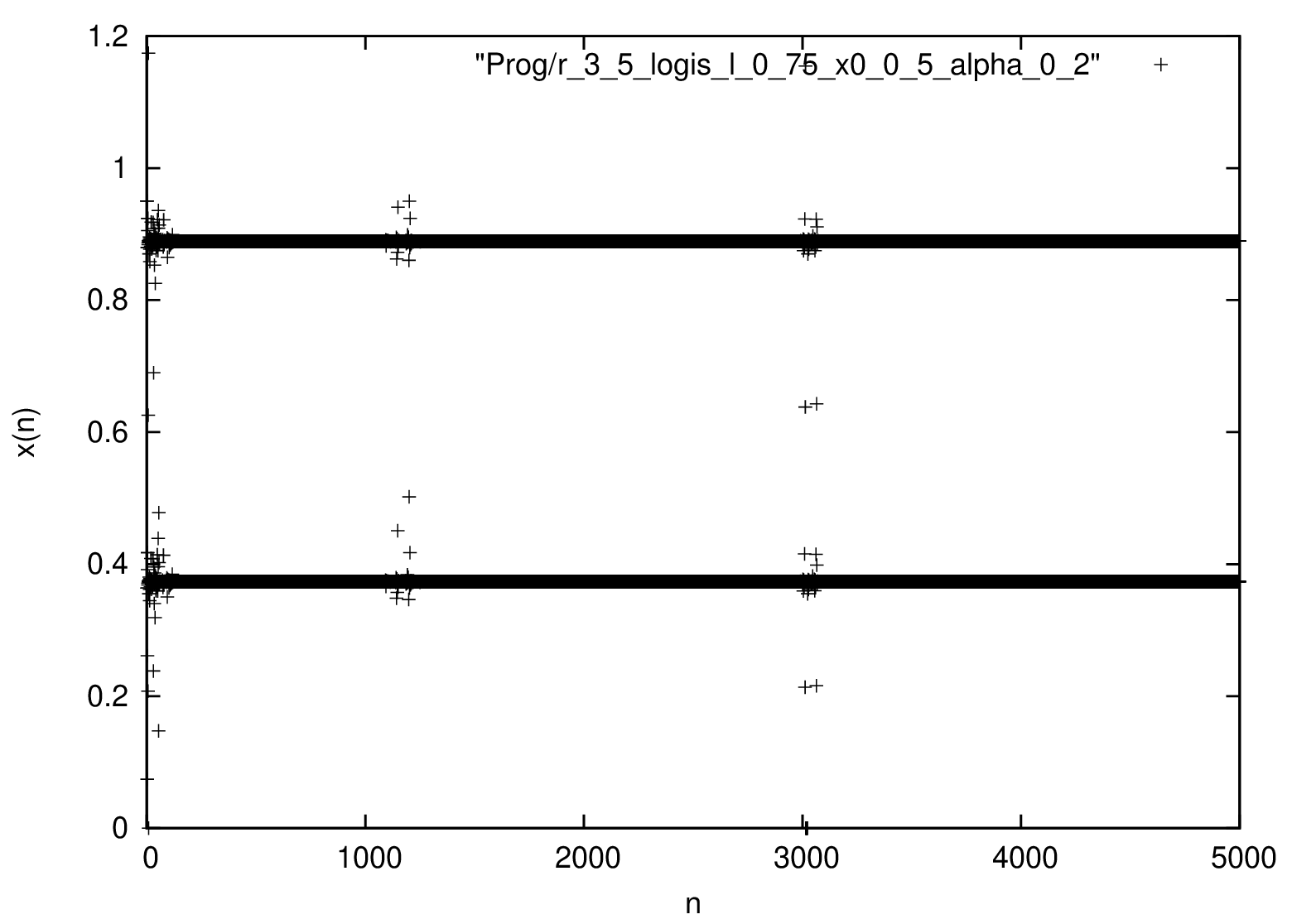}
\includegraphics[height=.125\textheight]{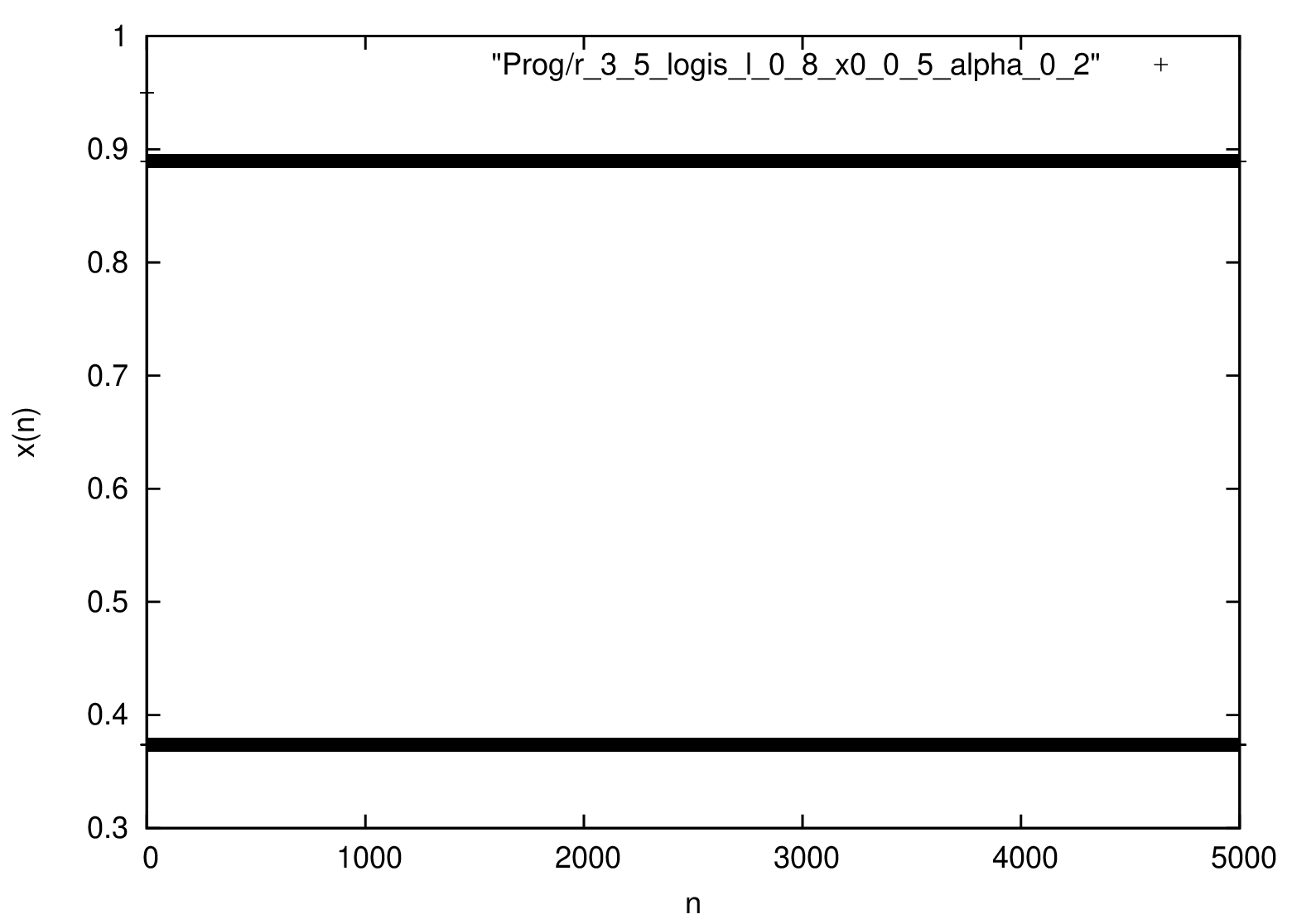}
\caption{
Model \eqref{eq:TOCKk} with $f=f_2$ from \eqref{eq:logistic},  $\alpha=0.2$, $m=1$, $d=2$, $x_0=0.5$ and (from left to right) 
$l=0.75,0.8$. 
}
\label{figure9new}
\end{figure}

\begin{example}
\label{ex:TOClogistic}
Consider stochastic TOC \eqref{eq:TOCKk} applied at alternate steps ($k=2$) 
to a logistic map $f_2$ satisfying \eqref{eq:logistic} with $r=3.5$. 
We can globally stabilize an unstable 2-cycle. 
Fig.~\ref{figure7new} shows convergence with  $\alpha=0$, $r=3.5$, $l=1.2$. For significatly smaller $l$, there is no convergence, and 
the effective stabilizing range for 
$l$ is narrow. Increasing to $\alpha=0.2$ in Fig.~\ref{figure9new} leads to a higher convergence speed, see 
Fig.~\ref{figure9new} for fast stabilization of a 2-cycle for  $l=0.75$ and $l=0.8$. 
\end{example}

Now we proceed to PBC method. Examples \ref {ex:PBCBevH} and \ref {ex:PBCRicker}  illustrate Theorem \ref{thm:PBCd}, $m=1$, with $d=1$ and $d=2$, 
respectively.

\begin{example}
\label{ex:PBCBevH}
Consider stochastic PBC \eqref{eq:kPBC} applied at every step to a Maynard-Smith model $f_3$ satisfying \eqref{eq:BevHolt} 
and note that this model is chaotic for $\alpha=0$,  $l=0$.  
There are two positive equilibria at $x =  2$ and $x = 4$, and $f_3'(2)=\frac 73 > 1$.  
Note also that even local stabilization is not possible for any $\alpha\in (0, 1)$ in the absence of noise ($l=0$).  
Fig.~\ref{figure11new} illustrates local stabilization of the equilibrium $x\equiv 2$  with $\alpha=0.8$ 
and  $l=0.4, 0.5, 0.7$. 
The solution $x\equiv 4$ is stable for $l=0$, and there is no stabilization of $x\equiv 2$ for $l=0.4$, we observe wandering 
between the two equilibria for $l=0.5$ and stabilization of $x\equiv 2$ for $l=0.7$.

\begin{figure}[ht]
\centering
\includegraphics[height=.125\textheight]{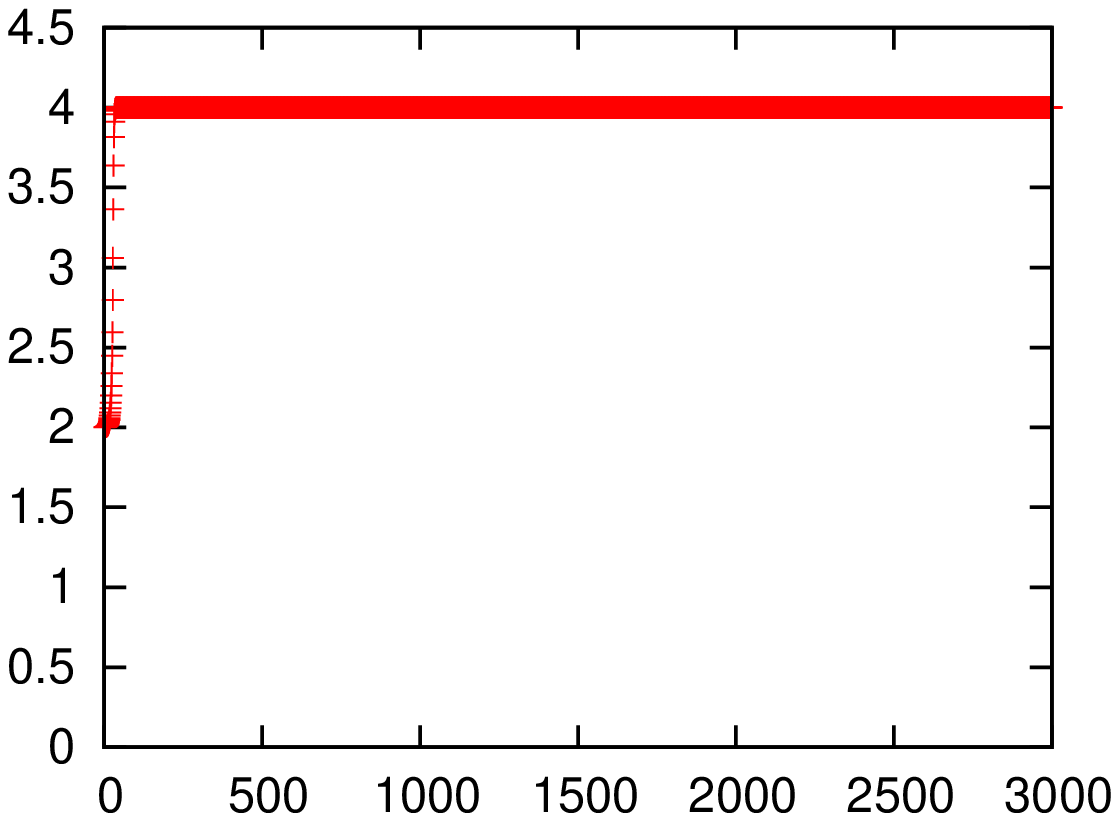}
\includegraphics[height=.125\textheight]{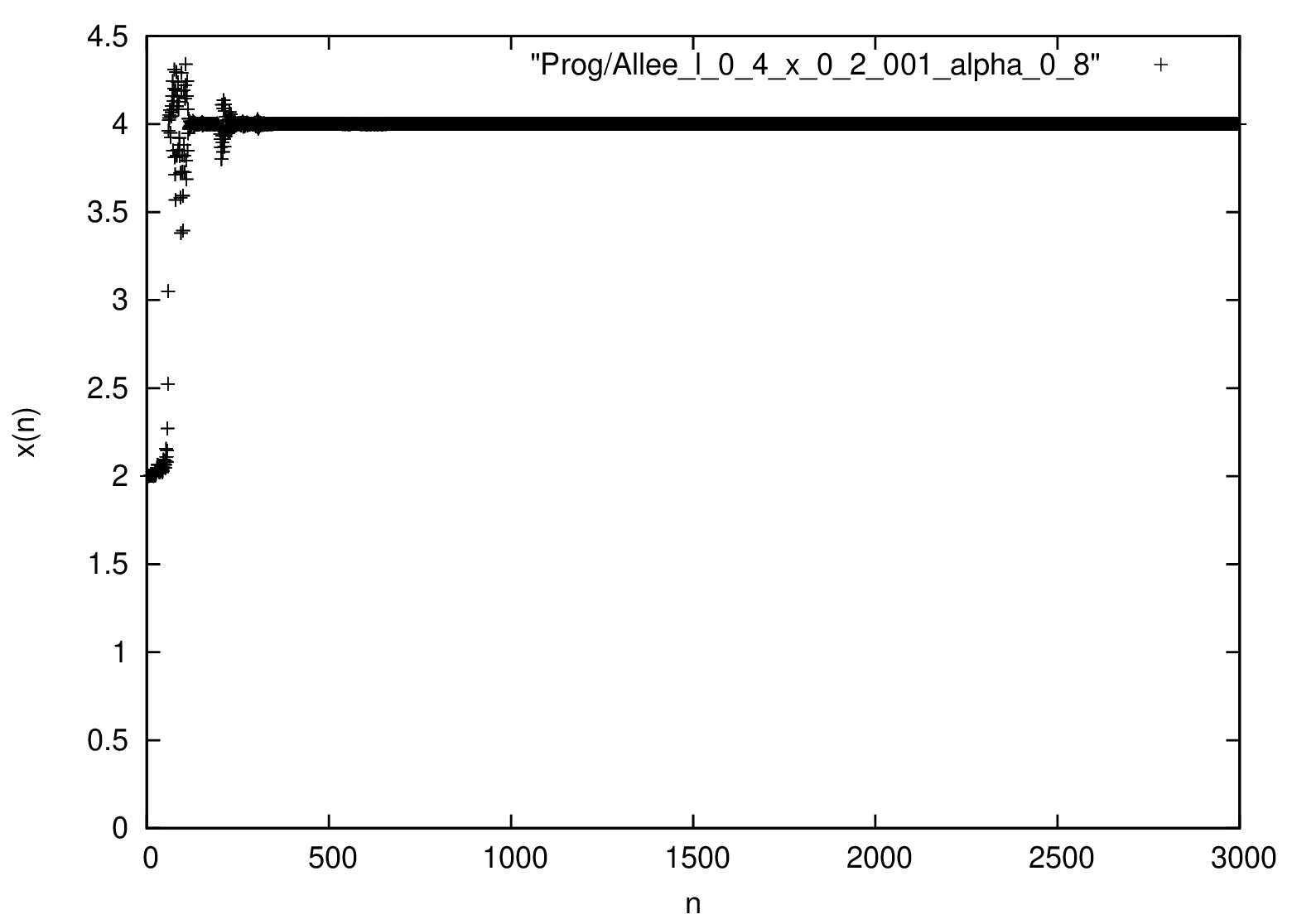} \\
\includegraphics[height=.125\textheight]{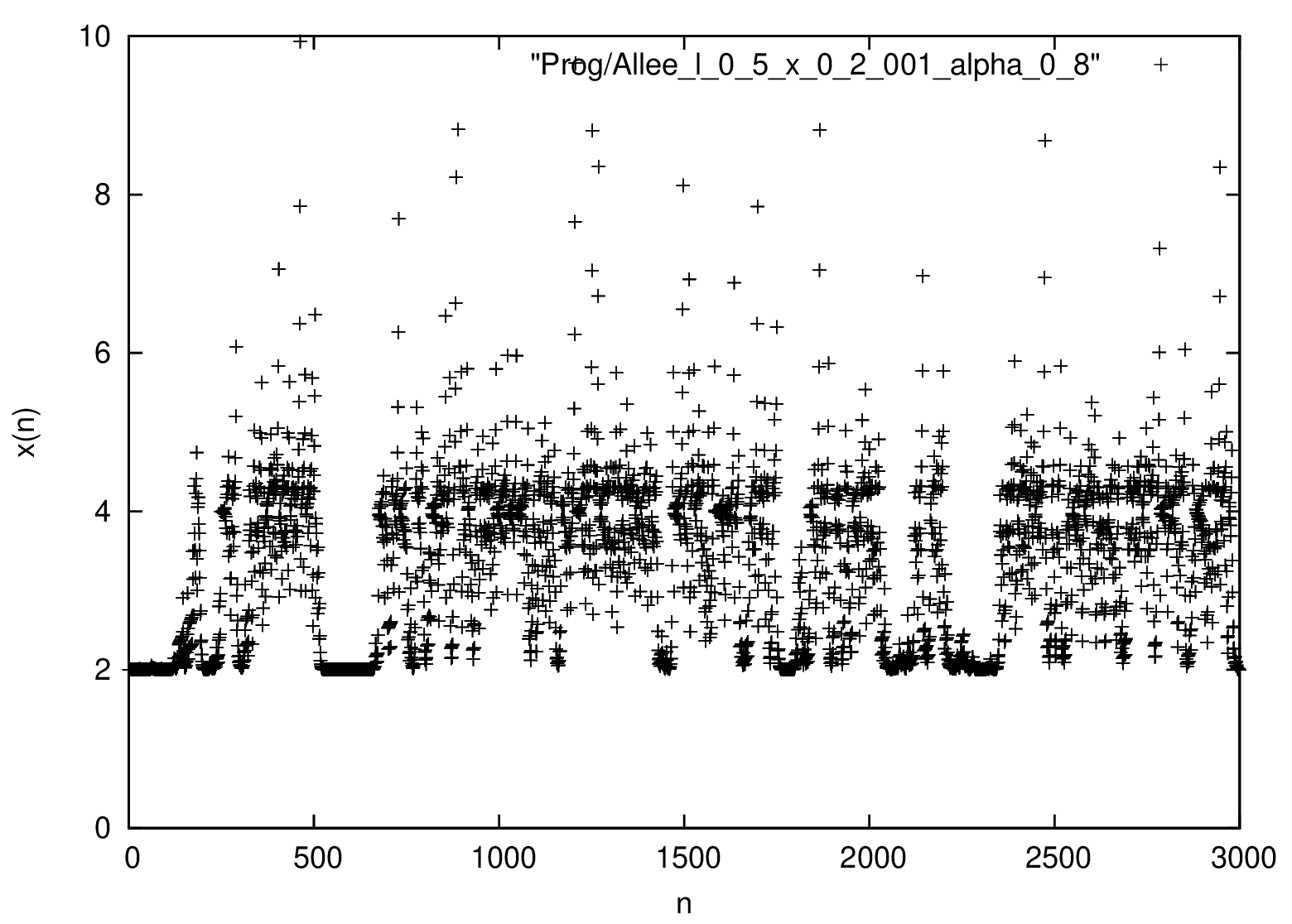}
\includegraphics[height=.125\textheight]{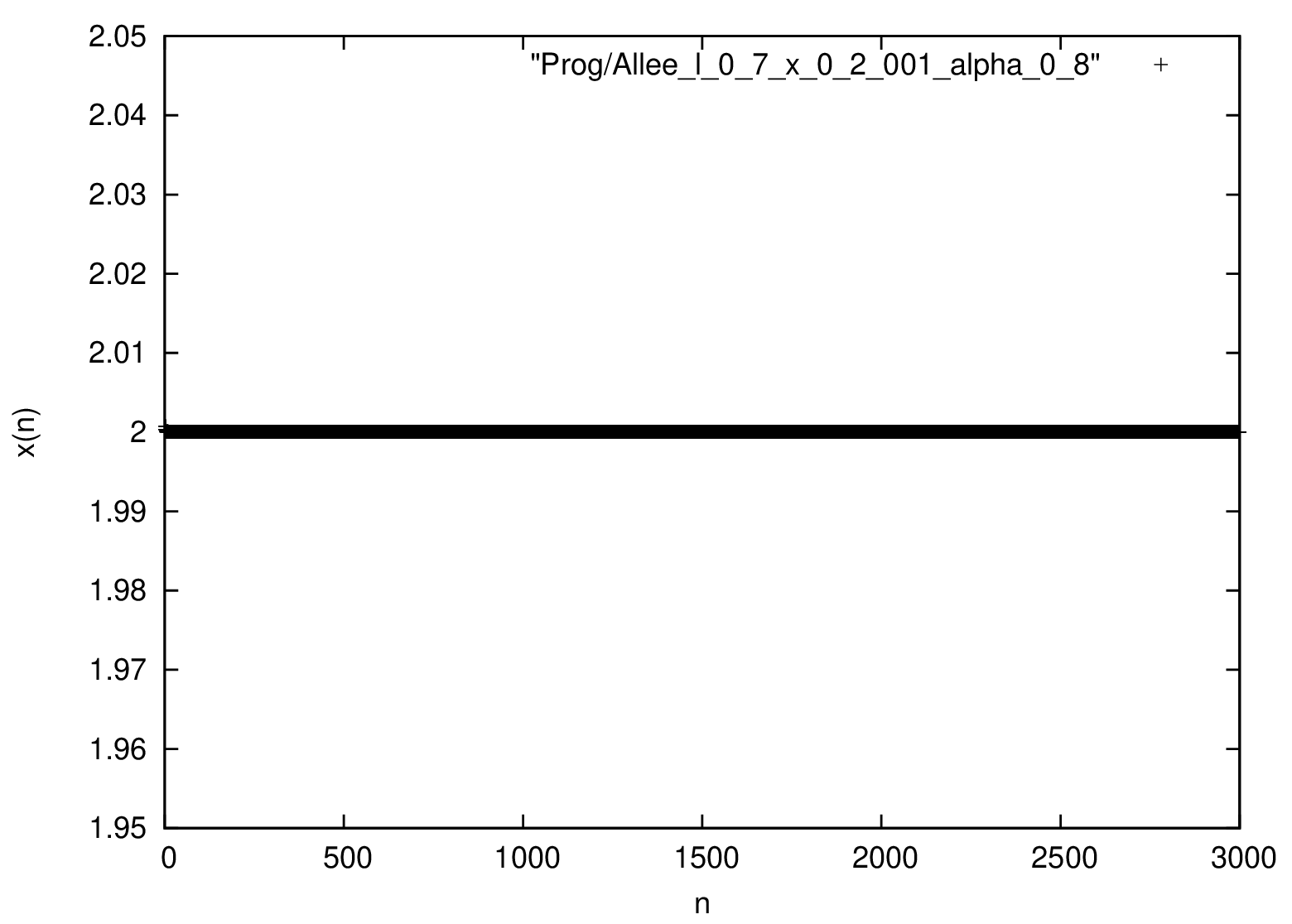}
\caption{Model \eqref{eq:kPBC} with $f=f_3$ from \eqref{eq:BevHolt}, $k=1$, $\alpha=0.2$, $x_0=2.001$
and (top) $l=04$, $l=0.4$ and (bottom) $l=0.5$, $l=0.7$.}
\label{figure11new}
\end{figure}

\end{example}

\begin{example}
\label{ex:PBCRicker}
Let us illustrate Theorem \ref{thm:PBCd} for $d=2$, $m=1$.
Applying pulsed stochastic PBC to stabilize a 2-cycle, we consider \eqref{eq:PBCmd},  
``delayed'' stabilization, applied to a Ricker map $f_1$ satisfying \eqref{eq:ricker} with $r=3.2$, $\alpha=0.4$. 
Fig.~\ref{figure14new} shows how an appropriately chosen noise intensity $l$ leads to global 
stabilization of the $2$-cycle $\{K_1\approx 0.11,K_2\approx 1.89\}$, and $l$ changing
 from zero (no noise) to $l=0.45$.

\begin{figure}[ht]
\centering
\includegraphics[height=.125\textheight]{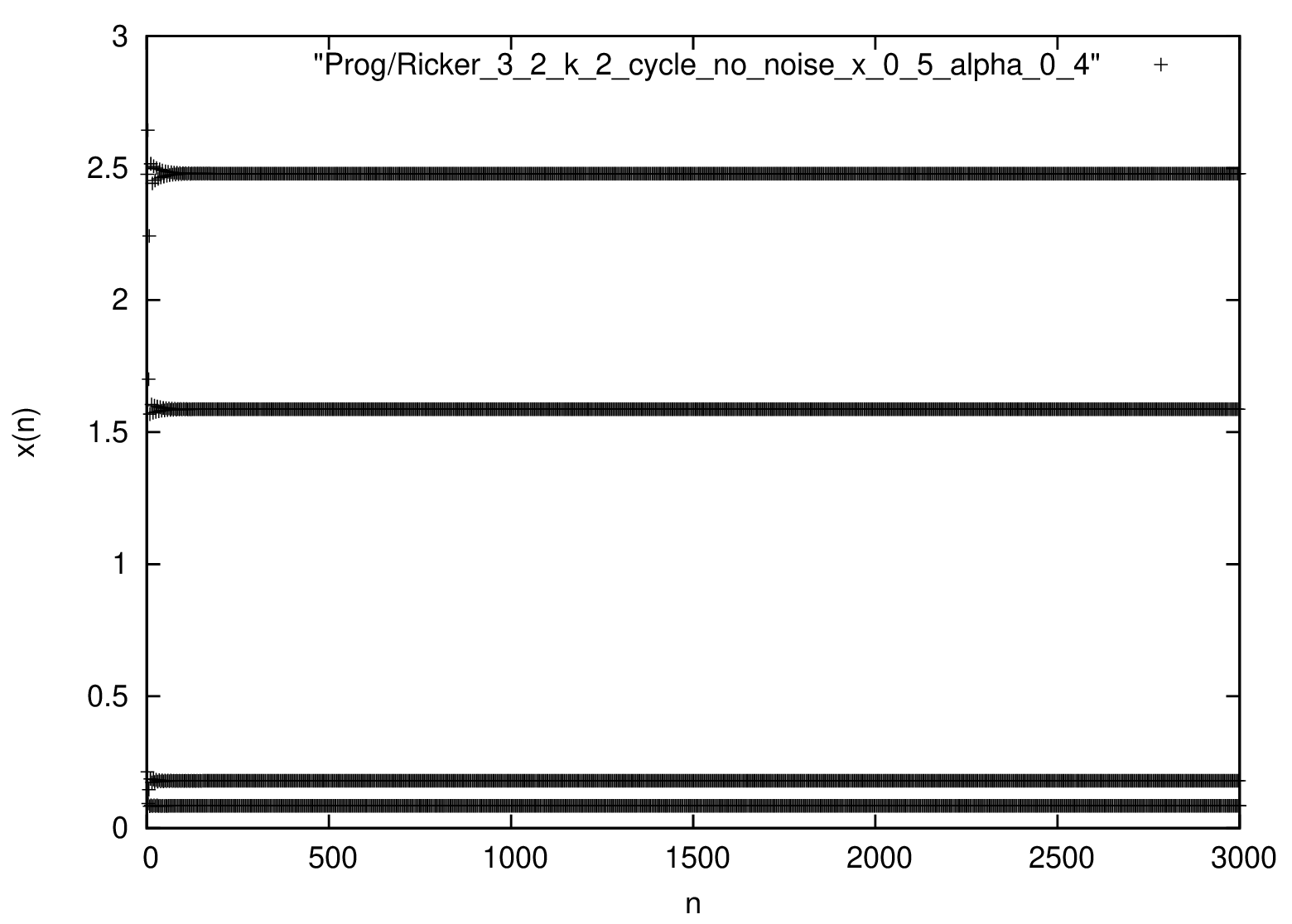}
\includegraphics[height=.125\textheight]{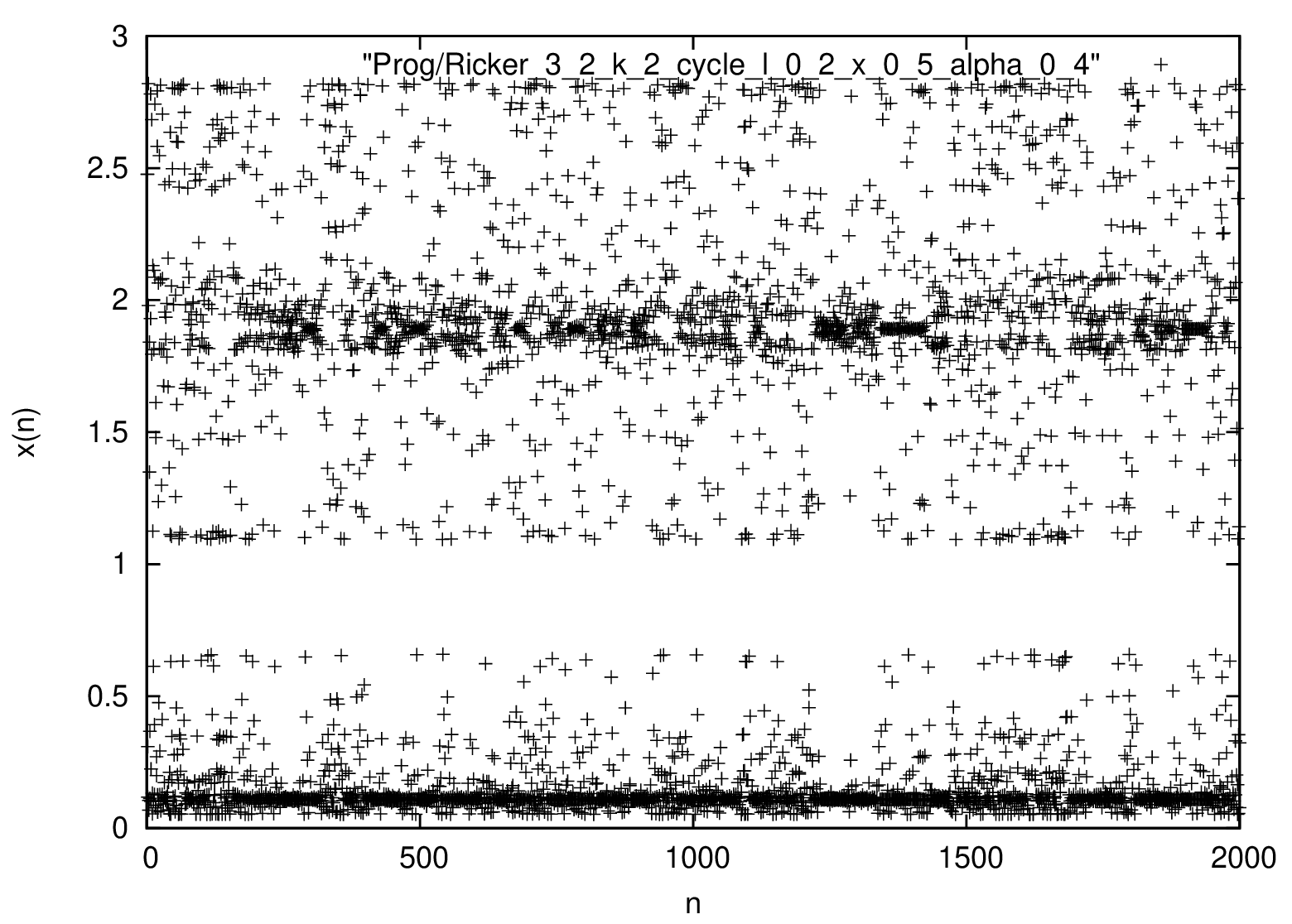} \\
\includegraphics[height=.125\textheight]{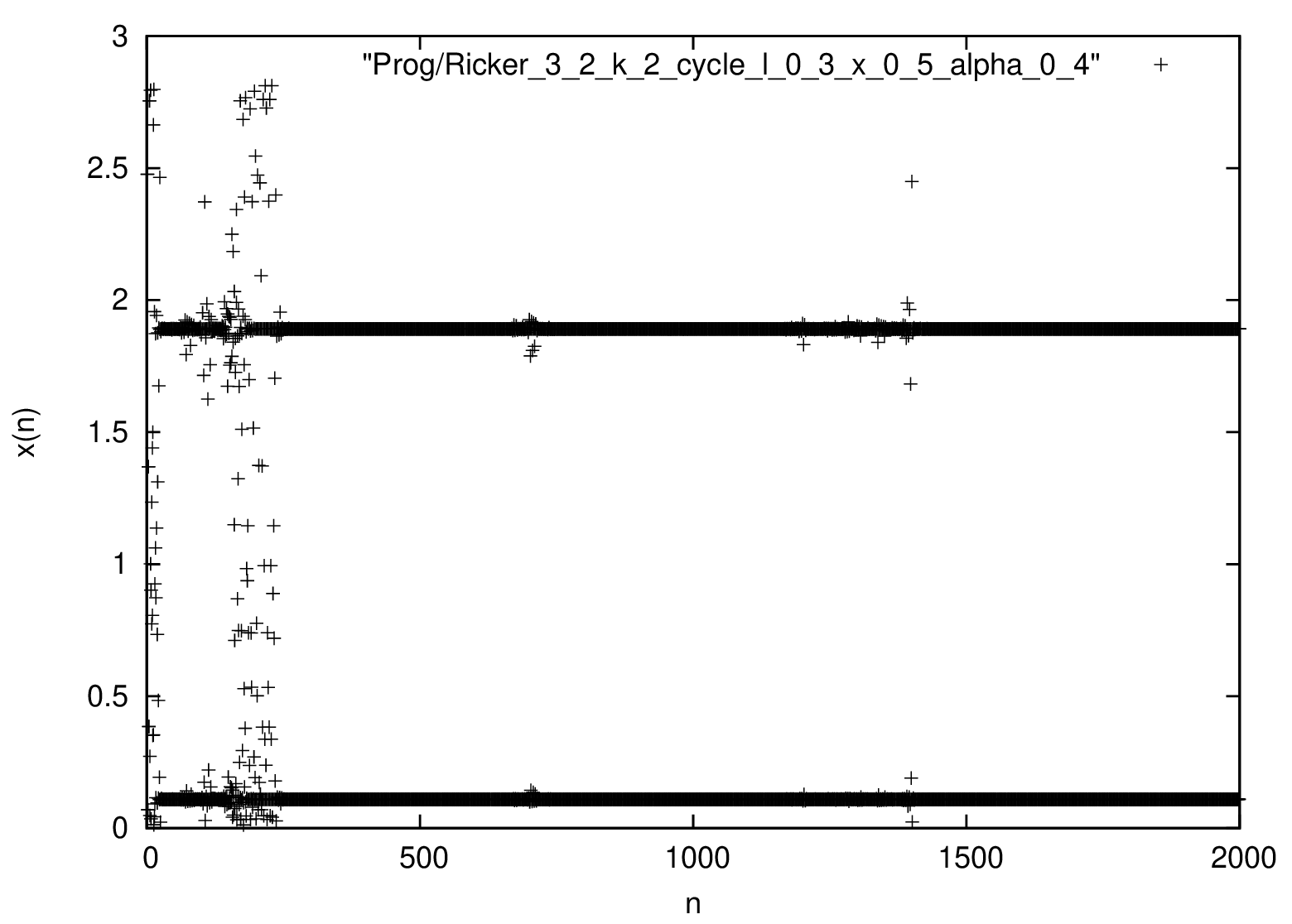}
\includegraphics[height=.125\textheight]{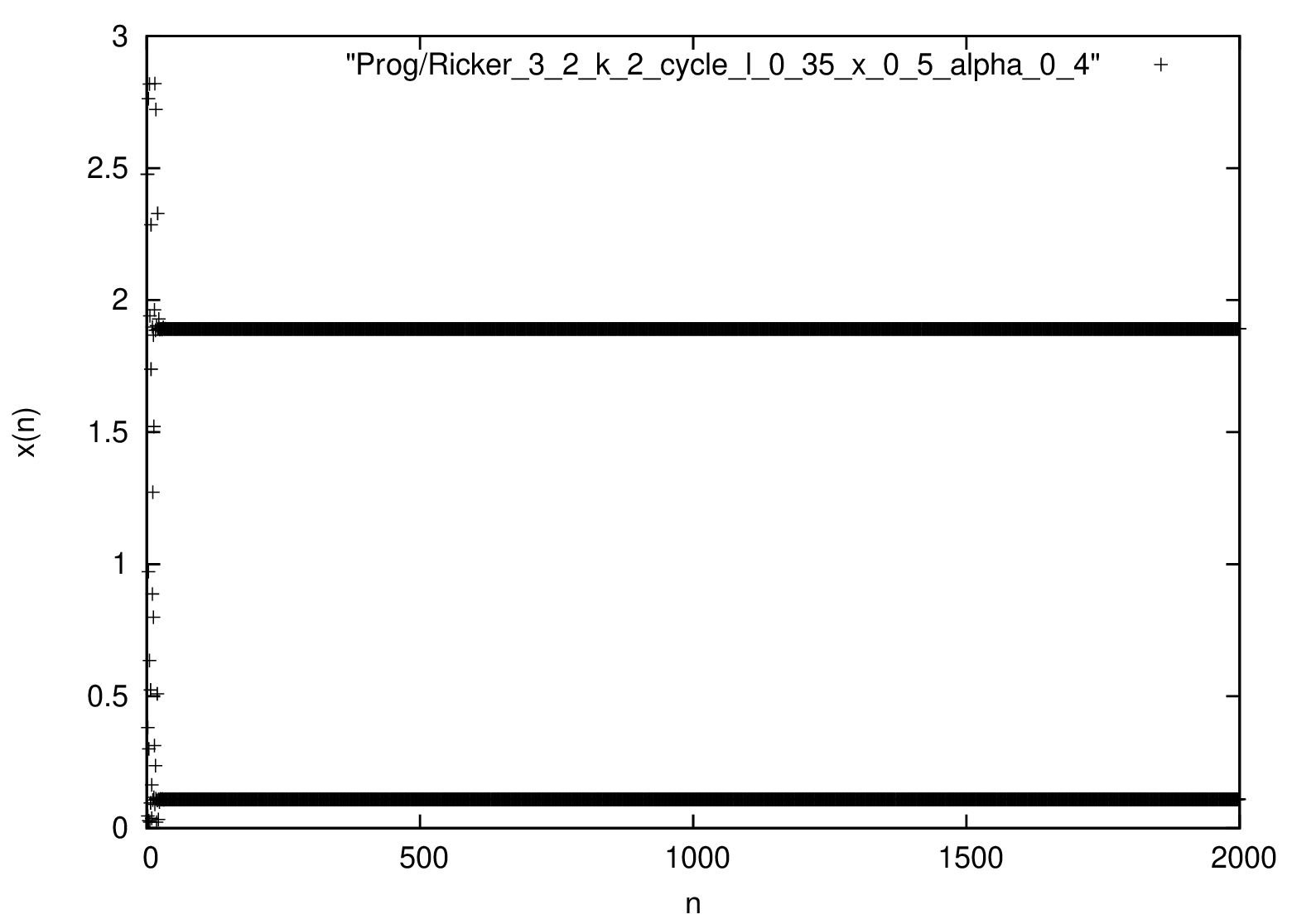} \\
\includegraphics[height=.125\textheight]{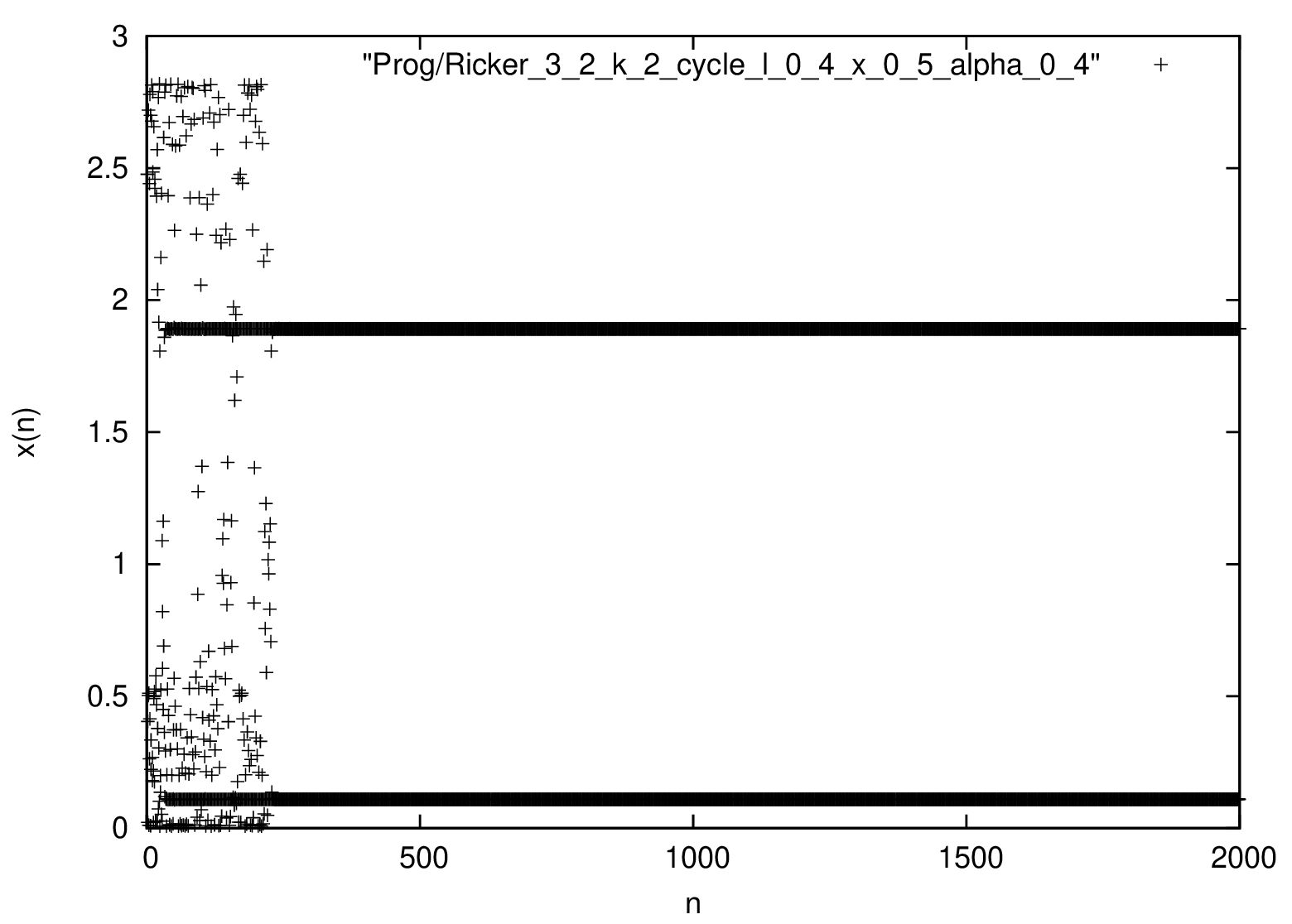}
\includegraphics[height=.125\textheight]{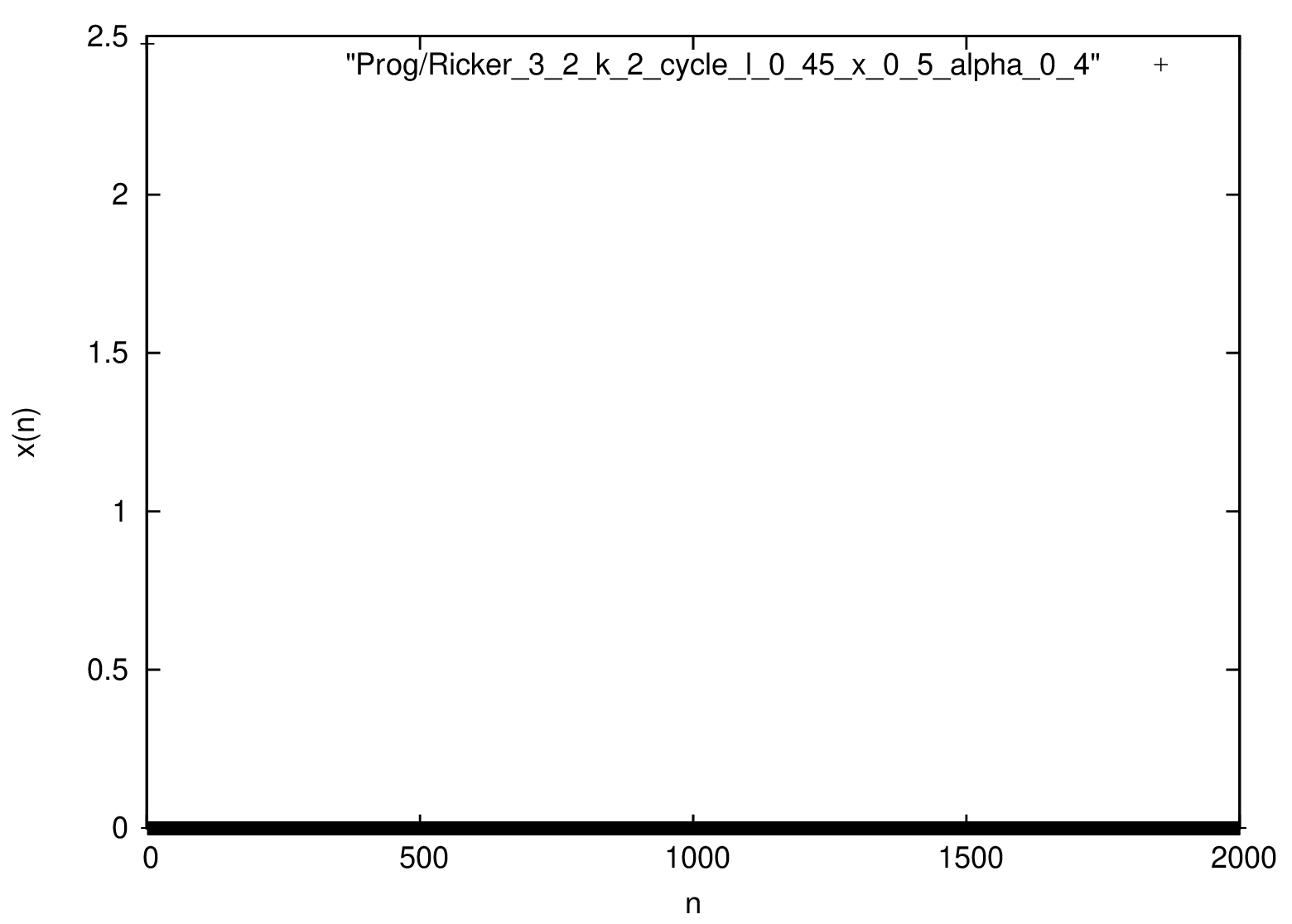}
\caption{Model \eqref{eq:PBCmd}  with  $f=f_1$ from \eqref{eq:ricker}
with $r=3.2$, $k=2$, $\alpha=0.4$, $x_0=0.5$ and (from left to right, top to bottom) $l=0,0.2,0.3,0.35,0.4,0.45$.
}
\label{figure14new}
\end{figure}

\end{example}

The next example illustrates Theorem~\ref{thm:fK<>k}.

\begin{example}
\label{ex:max}
Consider stochastic PBC \eqref{eq:kPBC} applied at each step to a Ricker map  
$f_1$ satisfying \eqref{eq:ricker} with $r=2.41$. The global Lipschitz 
constant is $L=1.5$. 
Note that according to Remark~\ref{rem:Schwarzian}, for $r=2.41$ we get
$$
f_{\alpha}^{\prime} (1) = -1.41(1-\alpha) + \alpha = 2.41 \alpha -1.41 > -1 $$ $$\Leftrightarrow   \alpha > \frac{0.41}{2.41} \approx  0.17012448.
$$
Thus the stabilization bound is $\alpha^*= 0.1701245$, and by \cite{BKR2016}, stabilization is achieved once $\alpha - l> \alpha^*$, $\alpha+l<1$. 
For $\alpha=0.3$, $l=0.24$, the first inequality is not satisfied $\alpha - l=0.06 < \alpha^*$, so our previous result in \cite{BKR2016} does not 
allow us to establish stability of the controlled model. It is possible to illustrate global stabilization of the equilibrium $K=1$  with $k=1$, 
$\alpha=0.3$ and  $l=0.24$.

We can also make $K=1$ stable with a stochastic pulsed control.
Global stabilization by pulsed stochastic PBC \eqref{eq:kPBC} applied at alternate steps ($k=2$) 
to a Ricker map $f_1$ satisfying 
\eqref{eq:ricker} with $r=2.2$ is demonstrated in Fig.~\ref {figure16new}. 
Only local stability conditions from Theorem~\ref{thm:mathcalA} hold, but global stability is observed. 
This indicates a possible direction for future research.

\begin{figure}[ht]
\centering
\includegraphics[height=.125\textheight]{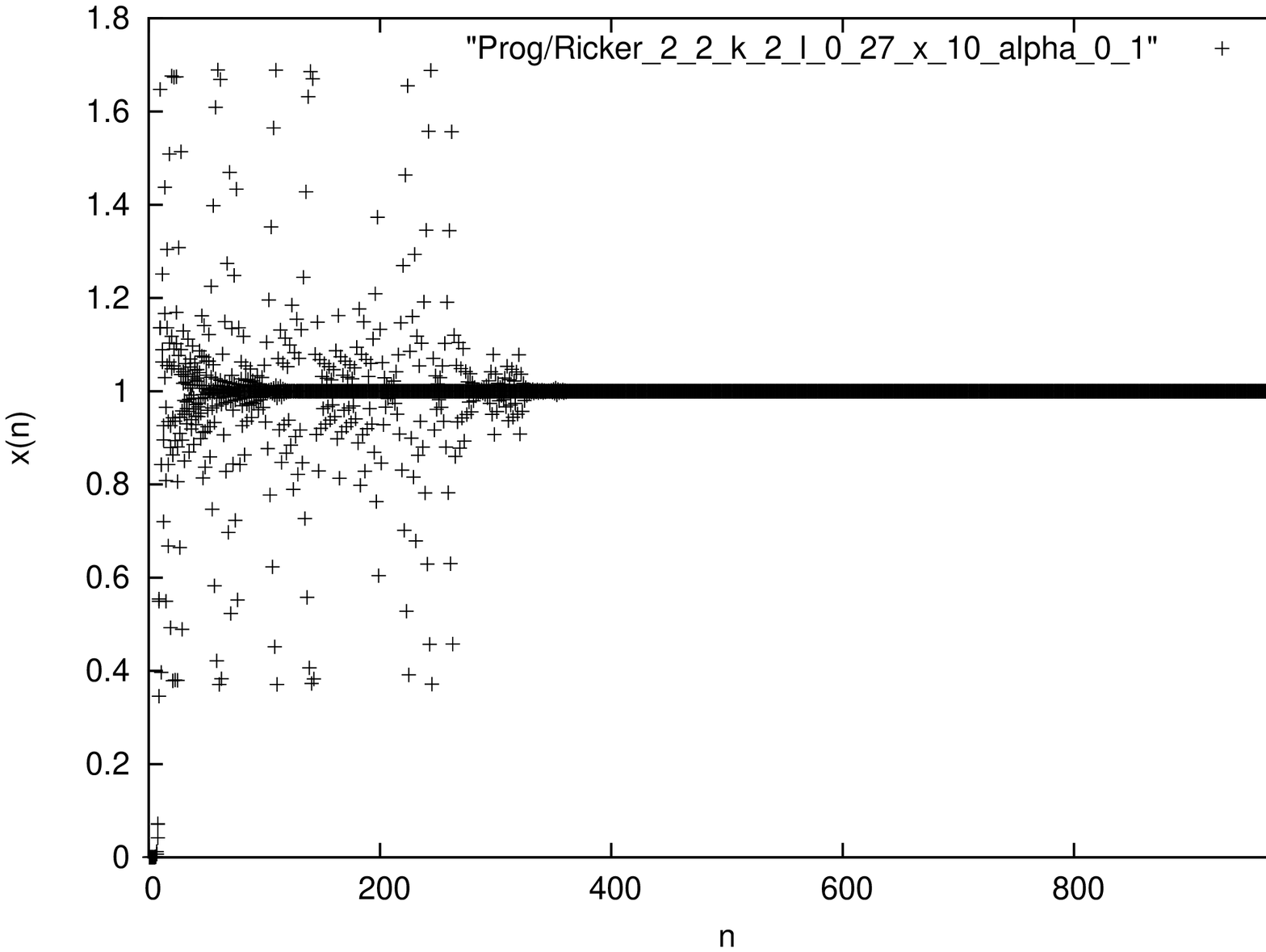}
\includegraphics[height=.125\textheight]{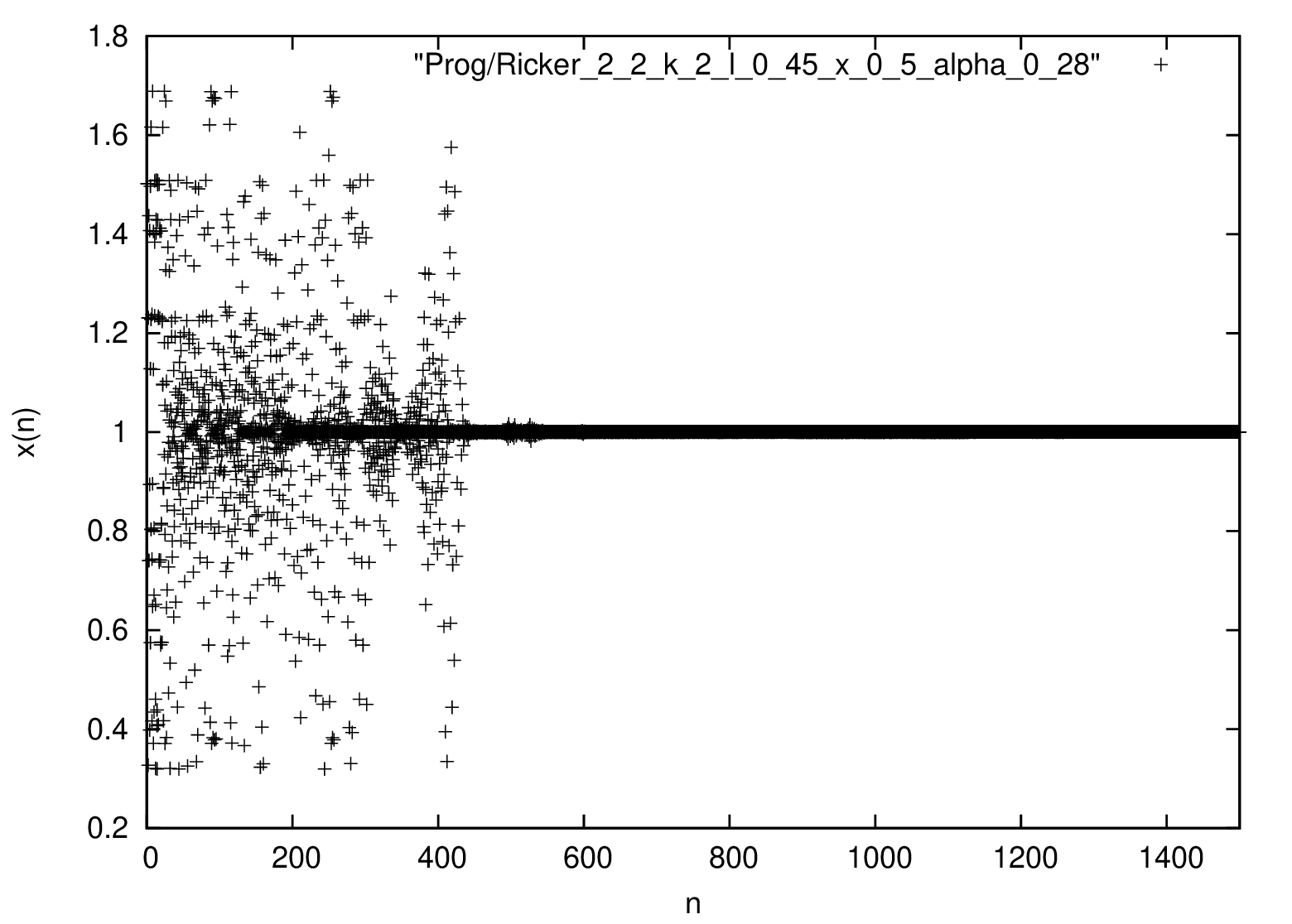}
\\
\includegraphics[height=.125\textheight]{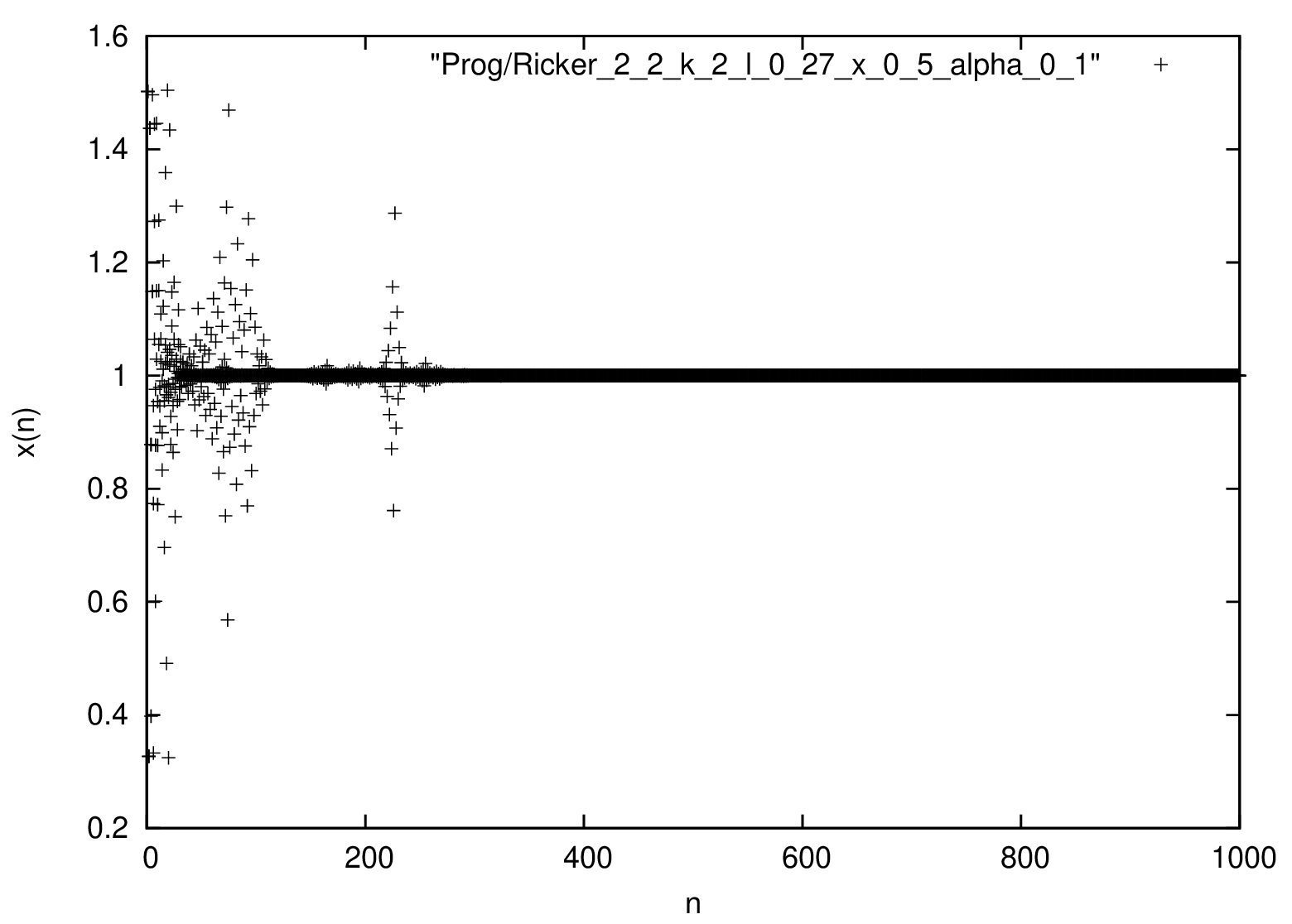}
\includegraphics[height=.125\textheight]{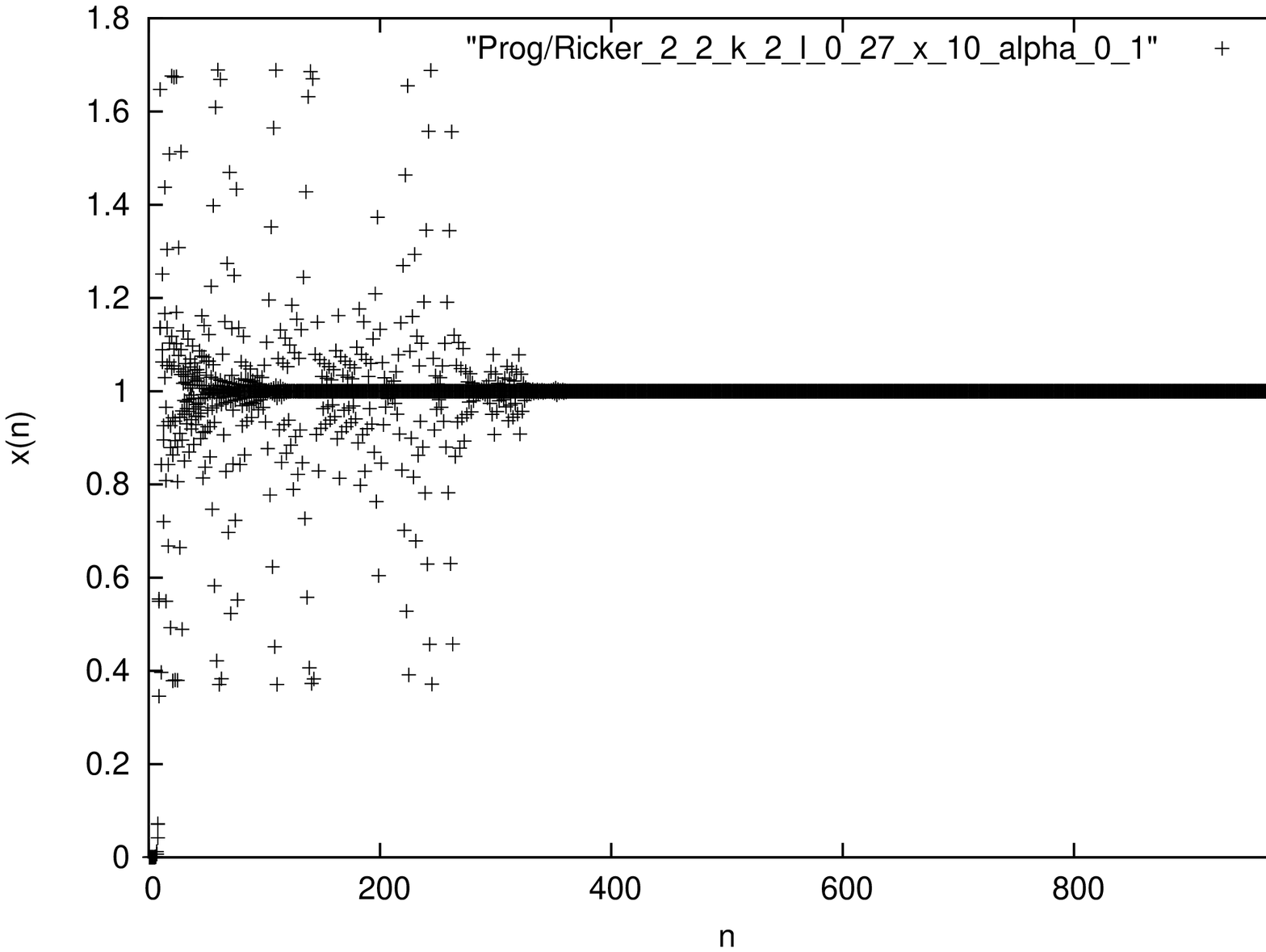}
\caption{Model \eqref{eq:kPBC}  with  $f=f_1$ from \eqref{eq:ricker}, $k=2$, $r=2.2$.
and (top left) $\alpha=0.28$,
$l=0.27$, $x_0=0.5$, (right) $l=0.45$, $x_0=0.5$,
(bottom) $\alpha=0.1$, $l=0.27$, (left)
$x_0=0.5$, (right) $x_0=10$.
}
\label{figure16new}
\end{figure}
\end{example}

\section[4]{Conclusions}
\label{sec:summary}
\subsection{Summary of results}

We have presented a general framework for stabilization of iterative systems by the application of control. 
The method demonstrates how the effective range of control parameters can be extended by stochastic perturbation.

Our results may be summarized as follows:
\begin{itemize}
\item
We present a broad general characterization of control 
that allows us to consider mechanisms incorporating both deterministic and stochastic components. 
This characterization includes stochastic forms of Prediction-Based Control and Target-Oriented Control; 
the latter is considered here for the first time.
\item
Both regular (applied at each step) and pulsed (applied every $k$th step) 
types of control were investigated. 
Pulsed TOC applied to one-dimensional models has not been studied before. 
Further, while pulsed stabilization of an equilibrium or a cycle 
by deterministic control has attracted some attention (see for example \cite{BL2012,LizPotsche14}), 
the use of stochastic PBC and TOC for $k$-cycle stabilization is novel.
\item
The analysis extends that of \cite{BKR2016}, which also showed how the effective 
range of stabilizing control parameters may be extended by the introduction of noise. 
However, in this article we additionally explore pulsed stochastic control and the stabilization 
of unstable $k$-cycles.
\end{itemize}

\subsection{Future directions for research}

Future research will follow one (or more) of the four directions:
\begin{enumerate}
\item
{\bf Study of sharp and/or global stabilization conditions.} 
Note that some of our results are essentially local, and we observe this in simulation. 
Nevertheless, we believe it is  still possible to get global stabilization results for the PBC method 
if we introduce some additional restrictions on $f$ and choose the noise intensity appropriately.
Example~\ref{ex:max} in fact illustrates the global stabilization for a wide range of parameters,
however theoretical justification is yet to be obtained.
Also, the conditions for stabilization presented in this article are sufficient but not necessary. 
It is desirable to obtain necessary 
and sufficient conditions, as are available in the deterministic setting: see \cite{TPC,Liz2007,FL2010}. 
\item
{\bf Explore the dependency of parameter bounds on the noise distribution type.} 
Most of our examples assume Bernoulli-distributed noise perturbation. 
Creating a library of sufficient estimates for control parameters and noise 
amplitudes under various types of noise distribution would be interesting and useful. 
For stabilization with noise only, such results can be found in \cite{BR2020}.
\item
{\bf Generalization of scalar results to systems or higher-order difference equations.}
In \cite{BF2015} and \cite{BF2017}, stabilization of high-order and vector difference equations was considered. 
This allowed to analyse stage-structured and delayed population dynamics models.
Even without formal control, the introduction of noise can improve population dynamics, for example, reduce oscillation amplitudes.
The next stage of research is to 
incorporate the ideas and methods of the present paper in the controlled systems of difference equations.  
For example, if a system describes a network, together with stabilization of periodic orbits \cite{LizPotsche14}, 
the problem of synchronization is 
of importance, especially if there are delays and/or stochastic component in communications \cite{Nag2016,Porfiri_SIAD}. 
\item
{\bf Qualitative analysis of positive effect of noise: ecological perspective.} 
Finally, it would be interesting, similarly to the present paper, to qualitatively evaluate possible positive effect 
of noise on stability and population survival, see \cite{Schreiber} for analysis of population interactions and 
the influence of stochasticity on survival, extinction, and coexistence.
\end{enumerate}


\section*{Acknowledgments}

The authors are grateful to Prof. Gregory Berkolaiko for the fruitful discussion of application of the Law of Large Numbers to PBC method and to two
anonymous reviewers whose thoughtful and valuable comments significantly contributed to the present form of the paper.   
The first author was  partially supported by the NSERC research grant RGPIN-2020-03934.

\section*{Data Availability Statement} 
The data that support the findings of this study are available from the corresponding author upon reasonable request.



\appendix
\label{sec:ap}

\setcounter{equation}{0}

\renewcommand{\theequation}{A.\arabic{equation}}


\subsection{Proof of Theorem \ref{thm:0equil}}
Consider two cases in turn: (a) $M L^{k-1}<1$ and (b) $M L^{k-1}\ge 1$. 
\bigskip

\noindent \underline{Case (a):} If $M L^{k-1}<1$ then \eqref{cond:EGLitsh} holds for $u_1=u_0$. Since Eq.  \eqref{cond:lipschg1} in Assumption \ref{as:g} holds with $L\ge 1$, we must have $M<1$. Fix $u\in (0, u_0]$ and set 
$\delta_0\le u/L^{k-1}.$ Then, for $|z_0|\le \delta_0$, 
\begin{eqnarray*}
|z_1|&\le& L|z_0|\le\frac {u}{L^{k-2}} \le u, \\
|z_2|&\le& L|z_1|\le L^2|z_0|<\frac {u}{L^{k-3}}\le u,\\
 &\vdots&\\
|z_{k-1}|&\le &L|z_{k-2}|\le L^{k-1}|z_0|\le u,\\
|z_{k}|&\le &M|z_{k-1}|\le ML^{k-1}|z_0|\le Mu<u,\\
|z_{2k}|&\le& \left( ML^{k-1} \right)^2 u\le u.
\end{eqnarray*}
Reasoning recursively, we conclude that, for $i \in {\mathbb N}$,
$$
|z_{ik}|\le \left( ML^{k-1} \right)^i u\le u, $$ $$  |z_{ik+j}|\le (ML^{k-1})^i u\le u, \, j=1, \dots, k-1,
$$
where $ML^{k-1}<1$, so, $\lim\limits_{i\to \infty}|z_{ik+j}(\omega)|=0$ for all $\omega\in\Omega$,  which concludes the proof of Part (a).
\bigskip

\noindent \underline{Case (b):} 
Let $\alpha$ and $l$ be chosen as in \eqref{cond:EGLitshu0},
and  $u_1$ be defined as in \eqref{cond:EGLitsh}.
Fix $u\le u_1$ and  $\gamma\in (0, 1)$.

Applying Lemma~\ref{thm:Kolm} and condition~\eqref{cond:EGLitsh} we conclude that, for
$\varepsilon:=1/2[\lambda -(k-1)\ln L]>0$, there exists a random $\mathcal N=\mathcal N(\gamma,
\lambda, u) $ such that
$$
\left|\frac 1n\sum_{i=0}^n \ln \mathcal L(\alpha, \xi_{ik},u)+\lambda\right|<\varepsilon,
\quad n\ge \mathcal N.
$$
Then there exist
 a nonrandom $N=N(\gamma, \lambda, u)$ and  $\Omega_\gamma\subset\Omega$ with
$ \mathbb P(\Omega_\gamma)>1-\gamma$, such that
\begin{equation}
\label{est:kolmL}
\begin{array}{ll}
\displaystyle & \prod_{i=0}^n \mathcal L(\alpha, l, \xi_{ik}(\omega), u) \\ < & \exp\left\{-\frac{n}{2}(\lambda +(k-1)\ln L)\right\}, ~
n\ge N,~\omega\in\Omega_\gamma. 
\end{array}
\end{equation}
Suppose that, when $N$ is chosen so that \eqref {est:kolmL} holds, $\delta_0$ satisfies
\begin{equation}
\label{def:delta}
\delta_0\le \frac {u} {(\bar ML^{k-1})^{N}}\,, \quad  \mbox{where} \quad \bar M:=\max\{1, M\}. \end{equation}
Since $\bar M\ge 1$, we have,  for all $i=1, \dots, k-1$, 
\begin{align*}
& \bar M^N L^{N(k-1)-i}=\bar M^N (L^{k-1})^{N-\frac i{k-1}} \\ \ge  & \bar M^{N-i/k-1} (L^{k-1})^{N-\frac i{k-1}}\ge 1.
\end{align*}
By   \eqref{cond:lipschg1}, \eqref{eq:genGalpha} and  \eqref{def:delta} we have, for $|z_0|\le \delta_0\le u$,
 \[
|z_{1}|=|g(z_0)| \le  \frac {L u} {(\bar ML^{k-1})^{N} }\le  \frac u{(\bar ML^{k-1})^{N-\frac 1{k-1}} }\le u,
\]
and inductively,  for all $i=0, \dots, k-1$, 
\[
 |z_{i}|\le  \frac {L^i u} {(\bar ML^{k-1})^{N} }\le \frac u{(\bar ML^{k-1})^{N-\frac i{k-1}} }\le u.
 \]
So
\begin{equation*}
\begin{split}
& |z_{k}|\le \mathcal L(\alpha, l, \xi_k, u)|z_{k-1}| \\ \le & \frac {\bar ML^{k-1}u}{(\bar ML^{k-1})^{N} }=\frac u 
{(\bar ML^{k-1})^{N-1}}\le u.
\end{split}
\end{equation*} 
Similarly, for  any $j< N$,
\begin{align*}
 |z_{jk}| & \le  \mathcal L(\alpha, l, \xi_{jk}, u)|z_{jk-1}| \\ & \le \mathcal L(\alpha, l, \xi_{jk}, u) 
L^{k-1}|z_{j(k-1)}|\\
& \le    (L^{k-1})^{j} \prod_{i=1}^j \mathcal L(\alpha, l, \xi_{ik}, u)|z_{0}| \\ & < (\bar ML^{k-1})^{j}|z_{0}|<\frac 
u {(\bar ML^{k-1})^{N-j}}\le u.
\end{align*}
Denoting
$
\bar \lambda:=[\lambda -(k-1)\ln L]/2,
$
and applying \eqref{est:kolmL}, we get, on $\Omega_\gamma$,
\begin{align*}
|z_{kN}| & \le L^{(k-1)N}|z_{0}|\prod_{i=1}^N \mathcal L(\alpha, l, \xi_{ik}, u)\\
& \le |z_{0}|e^{((k-1)N)\ln L}e^{- \frac {[\lambda +(k-1)\ln L] N}2}= e^{-\bar \lambda N}|z_{0}|<u,
\end{align*}
and then, for each $j=1, 2, \dots, k-1$, 
\begin{align*}
 |z_{kN+j}| & <   L^j|z_{kN}|< e^{-\bar \lambda N}L^j|z_{0}| \\ & < e^{-\bar \lambda N}\frac u {\bar M^NL^{(k-1)N-j}}<u.
\end{align*}
Similarly, for any  $n=km+j$, where $j=1, 2, \dots, k-1$, $m>N$, 
we get, for $\omega\in\Omega_\gamma$,
\begin{align*}
 |z_{n}(\omega)| & 
 <  L^j|z_{km}(\omega)|  <e^{-\bar \lambda km}L^j|z_{0}|\le e^{-\bar \lambda (n-j)}L^k|z_{0}| 
\\ & 
< e^{-\bar \lambda (n-j)}\frac u {\bar M^NL^{(k-1)N-k}}<ue^{-\bar \lambda (n-j)}<u, 
\end{align*}
which implies that $\lim\limits_{n\to \infty}z_n(\omega)=0$ when $\omega\in\Omega_\gamma$ and concludes the proof of Case (b).
\qed

{\bf Remark.} If $ML^{k-1}>1$ but $M<1$, we may use much bigger initial interval  than it was suggested in the proof of Theorem 
\ref{thm:0equil}
assuming $\delta_0\le \frac{u}{ML^{k-1}}$ for part (a) and $\delta_0\le \frac{u}{(ML^{k-1})^N}$ for part (b) of the proof. However, in this case we need  to change model \eqref{eq:genGalpha} slightly, considering instead
\begin{equation}
\label {eq:genGalphamod}
z_{n+1}=\left\{\begin{aligned} G(z_n, \alpha, l, \xi_{n+1}), &\quad n=k(s-1), \quad  s\in   \mathbb N;\\
g(z_n), &\quad \mbox{otherwise}, \quad z_0 \in \mathbb R.
\end{aligned}
\right.
\end{equation}
In \eqref {eq:genGalphamod}, the application of stochastic control starts from $n=0$ rather than $n=s-1$. 
Note that in both cases, the solution $z_n$ remains in $[-u, u]$ for all $n\in \mathbb N$. The proof differs  only in 
the estimation of the first $k$ (respectively $Nk$) iterations. 

\subsection{Proof of Lemma \ref{lem:Lipschfd}}

\label{subsubsec:proofLfd}
Let $i=1$, other cases are similar. Note that $\prod_{s=j}^d  \max\{L_{s}, 1\}\ge 1$ for $j=0, \dots, d$, and
$$
K_1=K_{d+1}=f(K_d)=
f^{d-j}(K_{d-j+1})=
f^{d}(K_1).
$$
For $ |x-K_1|\le u_0\left( \prod_{s=1}^d  \max\{L_{s}, 1\} \right)^{-1}$, we have
\begin{align*}
|f(x)-K_{2}| 
&\le \max\{L_1, \, 1\}|x-K_1| \\ & \le u_0\left( \prod_{s=2}^d  \max\{L_{s}, 1\}\right)^{-1}\le u_0,
\end{align*}
and, inductively, for each  $j=0, \dots, d-1$,
\begin{align*}
& |f^{d-j}(x)-K_{d-j+1}|\le \prod_{s=1}^{d-j}  \max\{L_s, \, 1\} |x-K_1|\\ \le & u_0\left( \prod_{s=d-j+1}^d  \max\{L_s, \, 1\}\right)^{-1}\le u_0,
\end{align*}
which, for $j=0$, implies \eqref{rel:Lipfd}.
\qed

\subsection{Proof of Lemma \ref{lem:expfd}}
Let $u(d)$ be defined as in \eqref{def:A(d)u(d)}.
Under Assumption \ref{as:fu0kcycle}, we have, for $x\in [K_i-u(d), K_i+u(d)]$, 
  \begin{align*}
 |f(x)-K_{i+1}|
 & \le |\mathcal A_{i}(x-K_i)|+ |\phi_{i}(x)|\\ 
& \le[ |\mathcal 
A_{i}|+ \psi_{i}(u)] |x-K_i|,
 \end{align*} 
 so \eqref{cond:fLip} holds for $L_i:= |\mathcal A_{i}|+ \psi_{i}(u)$. 
Acting as in the proof of Lemma \ref{lem:Lipschfd} we obtain that, for $j=1, \dots, d$, 
\begin{align*}
& |f^{j}(x)-K_{j+1}|\le \prod_{i=1}^{j}\left[|\mathcal A_{i}|+\psi_i(u)\right]|x-K_1|\\ \le &
\frac{u}{\prod_{i=j+1}^{d}  \max\{|\mathcal A_{i}|+\psi_i(u), \, 1\} }\le u,
\end{align*} 
where, notationally, $\prod_{i}^{j}\cdot=1 $ for any $i>j$. 
Now apply  \eqref {cond:fudcycle} recursively, for $|x-K_1|<u(d)$,  
\begin{align*}
& f^d(x)-K_1 \\ 
= & \mathcal A_{d}\biggl(\mathcal A_{d-1}(f^{d-2}(x)-K_{d-1}) \\ & +\phi_{d-1}(f^{d-2}(x))\biggr)+\phi_d(f^{d-1}(x))\\
\vdots&\\
= & \left(\prod_{i=1}^{d}\mathcal A_{i}\right)(x-K_1) +\sum_{j=1}^{d-1}\left( \prod_{i=j+1}^{d}\mathcal 
A_{i} \right) \phi_{j}(f^{j-1}(x)) \\ & +\phi_d(f^{d-1}(x)),
\end{align*}
where $f^0(x):=x.$ Define
\begin{equation*}
\label{def:phiuk}
\bar\phi(x):=\sum_{j=1}^{d}\left( \prod_{s=j+1}^{d}\mathcal A_{s}\right) \phi_{j}(f^{j-1}(x)). 
\end{equation*}
Acting as above  we get, for $j=1,\dots,d$,
\begin{align*}
& |\phi_{j}(f^{j-1}(x))| \\ \le & \psi_{j}( |f^{j-1}(x)-K_{j}|) |f^{j-1}(x)-K_{j}|\\
 \le  & \psi_{j}( |f^{j-1}(x)-K_{j}|) \left[ |\mathcal A_{j-1}|+\psi_{j-1}(u) \right] \\ & \times \left| f^{j-2}(x)-K_{j-1} \right|\\
 \le & \psi_{j}( |f^{j-1}(x)-K_{j}|)\prod^{j-1}_{i=1} \left[ |\mathcal A_{i}|+\psi_{i} (u) \right] |x-K_1|.
\end{align*}
So we can set 
\begin{align*}
\bar \psi(x):=b & \sum_{j=1}^{d}\biggl( \prod_{s=j+1}^{d}\mathcal A_{s}\biggr)\psi_{j}( |f^{j-1}(x)-K_{j}|)\\ & \times \prod^{j-1}_{i=1} \left(|\mathcal A_{i}|+\psi_{i} (u) 
\right),
\end{align*}
which completes the proof of  \eqref{cond:fudcycled}.
\qed

\subsection{Proof of Lemma \ref{lem:convmds+r}}
An application of Theorem \ref{thm:0equil} implies the existence of $s_0$ and $\Omega_\gamma\subset\Omega$ with $\mathbb P(\Omega_\gamma)\ge 1-\gamma$ such that, for sufficiently small $\delta>0$, $ s\ge s_0$, $\omega\in\Omega_\gamma$
$$
|x_{smd}(\omega)-K_1|<u_0L^{-m}(d)|1-\alpha+l|^{-m+1}.
$$ 
We need to show that, for each $s\ge s_0$, $j=qd+\bar j$, $\bar j=1, \dots, d-2$, and $q\le  m-1$,  we have $|x_{smd+j}(\omega)-K_{\bar j+1}|\le u_0$ when $\omega\in\Omega_\gamma$, which allows us to apply \eqref{cond:fLip} on each step.

Indeed, $|x_{smd+1}-K_{2}|=|f(x_{smd})-f(K_{1})|\le L_1|x_{smd}-K_1|\le u_0$. Reasoning inductively,  we have for each $j=1, \dots, d-1$ and $\omega\in\Omega_\gamma$,
\begin{equation*}
\begin{split}
&|x_{smd+j}(\omega)-K_{j+1}|=|f(x_{smd+j-1}(\omega))-f(K_{j})| \\ & \le \prod_{i=1}^{j}L_{i}|x_{smd}(\omega)-K_{1}|\le u_0,\\
&|x_{smd+d}(\omega)-K_{d+1}| \\ & \le  
|1-\alpha+l| L(d)|x_{smd}(\omega)-K_{1}|\le u_0,
\end{split}
\end{equation*}
and for $q=d+1,\dots,2d-1$ and again $\omega\in\Omega_\gamma$,
\begin{align*}
& |x_{smd+q}(\omega)-K_{q-d+1}| \\ 
\le & L_{q-d}|x_{smd+q-1}(\omega)-K_{q-d}| \\ \le & \dots \le \prod_{\theta=1}^{q-d}L_{\theta}|x_{smd+q-d}(\omega)-K_{1}|\\
\le & L(d) 
|1-\alpha+l|\prod_{\theta=1}^{q-d}L_{\theta}|x_{smd}(\omega)-K_{1}|\le u_0.
\end{align*}
Similarly, for $j=qd+\bar j$, $\bar j=1, \dots, d-2$, $q\le m-1$, and for $\omega\in\Omega_\gamma$,
\begin{align*}
& |x_{smd+j}(\omega)-K_{\bar j+1}| 
\\ \le & L_{\bar j}|x_{smd+j-1}(\omega)-K_{\bar j}|\\
\le & \dots  \le|1-\alpha+l|^q\prod_{\theta=1}^{\bar 
j}L_{\theta}L^q(d)|x_{smd}(\omega)-K_{1}|\\
\le & |1-\alpha+l|^{q}L^{q+1}(d)|x_{smd}(\omega)-K_{1}|   \\
\le & |1-\alpha+l|^{m-1}L^{m}(d)|x_{smd}-K_{1}|\le u_0.
\end{align*}
\qed

\subsection{Proof of Theorem \ref{thm:TOCglob}}
Consider \eqref{eq:TOCKk} with arbitrary $x_0$ and assume that \eqref{cond:ELalkTOC} holds.  
Then, for any $\gamma\in (0, 1)$, there exists $\Omega_\gamma\subset \Omega$ with $\mathbb P(\Omega_\gamma)>1-\gamma$ and $N\in\mathbb N$ such that for $\bar \lambda=\frac{1}{2} (\lambda-k\ln \bar L)$, we have on $\Omega_\gamma$,
  for $n=tk+j$, $j=1,2, \dots, k-1$, and $t=\lfloor n/k \rfloor>N$, where $\lfloor q \rfloor$ is an integer part of  $q\in [0,\infty)$,
  \begin{equation*}
  \begin{split}
 & \left( L^{k}|1-\alpha-l\xi_n| \right)^t\le   e^{-\bar \lambda t}, \\
&|x_{n}-K|\le L^{tk+j}|1-\alpha-l\xi_n|^t|x_0-K| \\
\le &  e^{-\bar \lambda t}L^k  |x_0-K|, 
\end{split}
\end{equation*}
which tends to zero as $n\to \infty$ (so that $k\to\infty$).
Here note that condition \eqref{cond:ELalkTOC} holds with $L$ substituted by the global constant $\bar L$,
which implies that $\lambda- k \ln \bar L>0$, and then $\bar \lambda>0$.

Analogously, consider  \eqref{eq:TOCmd} with arbitrary $x_0$, and assume that  \eqref{cond:ELalkTOCd} holds. 
Then
condition \eqref{cond:ELalkTOCd} is satisfied with $L(d)$
substituted by the global constant $\bar L(d)$, which implies that $\lambda- m \ln \bar L(d)>0$.
Applying Lemma \ref{lem:convmds+r} we get, 
for  $\bar \lambda=(\lambda-m\ln L(d))/2$, $n=tmd+j$,  $t=\lfloor u/(md) \rfloor>N$,  $j=qd+\bar j$, $\bar j=0, 1, \dots, d-1$, $q=0, 1, \dots, m-1$, 
\begin{equation*}
  \begin{split}
 & \left( L^{m}(d)|1-\alpha-l\xi_n| \right)^t\le   e^{-\bar \lambda t}, \\
&|x_{n}-K_{\bar j+1}|=|x_{tmd+j}-K_{\bar j+1}| \\ \le & |1-\alpha+l|^{m-1}L^{m}(d)|x_{tmd}-K_{1}|\\
\le & e^{-\bar \lambda t}   |1-\alpha+l|^{m-1}L^{m}(d)|x_0-K|\to 0, ~~ \mbox{as} ~~ m\to \infty,
\end{split}
\end{equation*}
and therefore $n \to \infty$.

If $ML^k<1$ (respectively, $ML(d)<1$), each of  limits above holds for all $\omega\in\Omega$.
\qed


\subsection{Proof of Theorem \ref{thm:fK<>k}}

Let $G(z, \alpha, l, v)$ be defined as in \eqref{def:PBCG} and denote, for simplicity, $G:=G(z, \alpha, l, v)$. Let $z\in [-u_0, 0]$ and $G>0$, then
\begin{align*}
| G| & =(1-\alpha-lv)g(z)+(\alpha+lv)z<(1-\alpha-lv)g(z)\\ & =(1-\alpha-lv) g(z)\le (1-\alpha-lv)L|z|.
\end{align*}
If $ G\le 0$, we have
$$
| G|=-(1-\alpha-lv)g(z)-(\alpha+lv) z<(\alpha+lv)|z|.
$$
Now, let $z\in [0, u_0]$ and $ G>0$, then
\[
| G|=(1-\alpha-lv)g(z)+(\alpha+lv)z<(\alpha+lv) |z|.
\]
If $ G\le 0$, 
\begin{align*}
| G|  & =-(1-\alpha-lv)g(z)-(\alpha+lv)z \\ & <(1-\alpha-lv)|g(z)|\le (1-\alpha-lv)L|z|.
\end{align*}
So \eqref{cond:GLitsh} holds  for ${\mathcal L}(\alpha, l, v, u)=\max_{|v|\le 1}\{(1-\alpha-lv)L, \,\alpha+lv \}$, and  \eqref{est:L} is satisfied for $M=\max\{(1-\alpha+l)L, \,\alpha+l \}$. 
Condition~\eqref{cond:EGLitshu0}  then takes the form
\[
\min\left\{-\mathbb E\ln |1+\alpha|, \,  -\mathbb E\ln |1-\alpha-l\xi|L\right\}\ge (k-1) \ln L.
\]
A direct application of  Theorem~\ref{thm:0equil} proves Part (ii).

The condition $ML^{k-1}<1$ takes the form $\max\{(1-\alpha+l)L, \,\alpha+l \}L^{k-1}<1$. In the case where $k=1$  and all $L>1$, this implies the conditions listed in \eqref{eq:tpia} in Part (i). In the case  where $k>1$ and for all $L$ satisfying $1<L^k<L+1$, this implies the conditions in \eqref{eq:tpib} in Part (i).  Another application of Theorem~\ref{thm:0equil} concludes the proof.
\qed


\begin{thebibliography}{99}


\bibitem{ABR2009} 
J.\,A.\,D.~Appleby, G. Berkolaiko and A.Rodkina, 
Non-exponential stability and decay rates in nonlinear stochastic difference equations 
with unbounded noise, {\em Stochastics} {\bf 81} (2009), 
99--127. 


\bibitem{Bashk}
I. Bashkirtseva, 
Mean-square analysis of stochastic cycles in nonlinear discrete-time systems with parametric noise, 
{\em J. Difference Equ. Appl.} {\bf 20} (2014), 
1178--1189.

\bibitem{Bashk1}
I. Bashkirtseva, E. Ekaterinchuk and L. Ryashko,  
Attractors of randomly forced logistic model with delay: stochastic sensitivity and noise-induced transitions,
{\em J. Difference Equ. Appl.} {\bf 22} (2016), 
376-–390.

\bibitem{Chaos2014}
E.~Braverman and B.~Chan,
\newblock Stabilization of prescribed values and periodic orbits with regular
and pulse target oriented control.
\newblock {\em Chaos}, 24(1):013119, 2014.

\bibitem{BF2015}
E. Braverman and D. Franco,  
Stabilization with target oriented control for
higher order difference equations.
\newblock {\em Physics Letters A} \textbf{379}(16) (2015), 1102--1109.

\bibitem{BF2017}
E. Braverman and D. Franco,
Stabilization of structured populations via vector
target-oriented control,
{\em Bull. Math. Biol.} {\bf 79} (2017), 1759-–1777.
	
\bibitem{BKR2016}
E. Braverman, C. Kelly and  A. Rodkina,
Stabilisation of difference equations with noisy prediction-based control,
{\it Physica D}  {\bf 326} (2016), 21--31.

\bibitem{BL2012}
E.~Braverman and E.~Liz,
\newblock On stabilization of equilibria using predictive control with and
without pulses.
\newblock {\em Comput. Math. Appl.} {\bf 64} (2012), 2192--2201. 

\bibitem{BR2019}
\newblock E. Braverman and A. Rodkina,
\newblock {Stochastic control stabilizing unstable or chaotic maps},
 {\em J. Difference Equ. Appl.} { \bf 25} (2019), 151-178.

\bibitem{BR2020}
E. Braverman and A. Rodkina,
Global stabilization and destabilization by the state dependent noise
with particular distributions,
{\em Physica D} {\bf 403} (2020), 132302,
https://doi.org/10.1016/j.physd.2019.132302.

\bibitem{Calatayud}
J. Calatayud,  J.-C. Cort\'{e}s, M. Jornet and L. Villafuerte, 
Random non-autonomous second order linear differential equations: mean square analytic solutions and their statistical properties, 
{\em Adv. Difference Equ.} 2018, Paper No. 392, 29 pp. 
 
\bibitem{Cortes}
J.-C. Cort\'{e}s, A. Navarro-Quiles, J.-V. Romero and M.-D. Rosell\'{o},  Full solution of random autonomous first-order linear systems of difference equations. Application to construct random phase portrait for planar systems, {\em Appl. Math. Lett.} {\bf  68} (2017), 150--156. 

\bibitem{Dattani}
J.~Dattani, J.\,C. Blake, and F.\,M. Hilker,
Target-oriented chaos control,
{\em Phys. Lett. A} {\bf 375}
(2011), 3986--3992.
 
\bibitem{Du}
N.\,H. Du and N.\,T. Dieu, 
Stochastic dynamic equations on time scales, {\em Acta Math. Vietnam.} {\bf 38} (2013), 
317-–338.

\bibitem{TPC}
D.~Franco and E.~Liz,  
A two-parameter method for chaos control and targeting in one-dimensional maps,
{\em  Internat. J. Bifur. Chaos Appl. Sci. Engrg. Int.} {\bf 23} 
(2013), 1350003, 11pp.

\bibitem{Franco2020}
D.~Franco , J. Per\'{a}n  and  J. Segura,
Stability for one-dimensional discrete dynamical systems revisited,
{\em Discrete Contin. Dyn. Syst. Ser. B}  {\bf 25} (2020), 
635--650.

\bibitem{FK}
H. Furstenberg and H. Kesten, 
Products of random matrices, {\em Ann. Math. Statist.} {\bf 31} (1960),
457--469.

\bibitem{Guzik}
G. Guzik,  
Asymptotic stability of discrete cocycles, {\em J. Difference Equ. Appl.} {\bf 21} (2015), 
1044-–1057.

\bibitem{Hasmin}
R.\,Z. Has'minski,
Stability of Systems of Differential Equations under Random Perturbations of Their 
Parameters, Nauka, Moscow, 1969, 367 pp.
 
\bibitem{Medvedev}
P. Hitczenko and G. Medvedev, Stability of equilibria of randomly perturbed maps, Discrete Contin.
Dyn. Syst. Ser. B 22 (2017), 269--281.

\bibitem{Kadiev}
R. Kadiev and P. Simonov,  
Initial data stability and admissibility of spaces for It\^{o} linear difference equations, {\em Math. Bohem.} {\bf 142} (2017), 
185--196.

\bibitem{Kapica}
P.~P.~Kapica,
Dynamical stability of pendulum with vibrating suspension,
{\em Journal of Experimental and Theoretical Physics} {\bf 21}
(1951),  588--597 (in Russian).

\bibitem{K} 
H. Kesten, Random difference equations and renewal theory for the
product of random matrices, {\em Acta Math.} {\bf 131} (1973), 207--248.


\bibitem{Liz2007}
E. Liz, 
Local stability implies global stability in some one-dimensional discrete single-species
models, {\em Discrete Contin. Dyn. Syst. Ser. B} {\bf 7} (2007), 191--199.

\bibitem{FL2010}
E.~Liz and D.~Franco, 
\newblock Global stabilization of fixed points using predictive control,
\newblock {\em Chaos} {\bf 20} (2010), 023124, 9 pages.

\bibitem{LizPotsche14}
E.~Liz and C.~P\"otzsche,
\newblock PBC-based pulse stabilization of periodic orbits,
\newblock {\em Physica D} {\bf 272} (2014), pp. 26--38.



\bibitem{Manjun}
G. Manjunath and H. Jaeger, 
The dynamics of random difference equations is remodeled by closed relations, 
{\em SIAM J. Math. Anal.} {\bf 46} (2014), 
459--483.

\bibitem{Nag2016}
M. Nag and S. Poria, 
Synchronization in a network of delay coupled maps with stochastically switching topologies, 
{\em Chaos Solitons Fractals} {\bf 91} (2016), 9--16.

\bibitem{Porfiri_SIAD}
M. Porfiri and I. Belykh, 
Memory matters in synchronization of stochastically coupled maps, 
{\em SIAM J. Appl. Dyn. Syst.} {\bf 16} (2017), 
1372--1396.

\bibitem{Rodina}
L.\,I. Rodina, On repelling cycles and chaotic solutions of difference equations with random parameters (Russian), {\em Tr. Inst. Mat. Mekh.} {\bf 22} (2016), 
227-–235.

\bibitem{Schreiber2001}
S.\,J. Schreiber, 
Chaos and population disappearances in simple ecological models, 
{\em J. Math. Biol.} {\bf 42} (2001), 
239--260. 

\bibitem{Schreiber}
S.\,J.~Schreiber, 
Coexistence in the face of uncertainty, {\em Recent progress and modern challenges in applied mathematics, modeling and computational science}, 349–-384, Fields Inst. Commun., {\bf 79}, Springer, New York, 2017.

\bibitem{Shaikhet}
L. Shaikhet,
Lyapunov Functionals and Stability of Stochastic Difference Equations, Springer-Verlag, London, 2011.


\bibitem{Shiryaev96} 
A.~N.~Shiryaev.
\newblock \emph{Probability}. 2nd edition,
\newblock Springer, Berlin, 1996.


\bibitem{Singer}
D. Singer, 
Stable orbits and bifurcation of maps of the interval, {\em SIAM J. Appl. Math.} {\bf 35}
(1978), 260--267.

\bibitem{Thieme}
H.R. Thieme, Mathematics in Population Biology, Princeton
University Press, Princeton, 2003.


\bibitem{uy99} 
T.~Ushio and S.~Yamamoto,   
\newblock Prediction-based control of chaos, 
\newblock {\em Phys. Lett. A} {\bf  264} (1999),  30--35.




\end{thebibliography}
\end{document}